\documentclass[11pt]{amsart}
\usepackage{latexsym,amsfonts,amssymb,multicol}
\usepackage{graphicx}
\usepackage[latin1]{inputenc}

\DeclareGraphicsExtensions{.jpg,.pdf}

\newtheorem{thm}{Theorem}
\newtheorem{prop}{Proposition}
\newtheorem{lem}{Lemma}
\newtheorem{cor}{Corollary}

\newtheorem{hp}{Assumption}

\renewcommand{\theequation}{\thesection.\arabic{equation}}

\newcommand{\supp}[1]{\supp(#1)}

\newcommand{\eop}{\hspace*{\fill} \ensuremath{\Box} \bigskip}

\def\Var{\mathop{\rm Var}\nolimits}%

\def \Sum{\displaystyle \sum}

\def \Sup{\displaystyle \sup}

\def\esp{\mathbb{E}}

\def\Diag{\mbox{\rm Diag}}

\def\ind{\mathbb{I}}

\def\Nit{\mathbb{N}}

\def\Rit{\mathbb{R}}

\def\1{\mathbb{I}}

\def\Prob{\mathbb{P}}
\def\prob{\mathbb{P}}

\def\supp{\mbox{\rm Supp}}

\renewcommand{\baselinestretch}{1.5}

\newcounter{numrems}

\def\Sum{\displaystyle \sum}

\def\esp{\mathbb{E}}
\def\Var{\hbox{\rm Var}}

\def\ind{\mathbb{I}}
\def\Nit{\mathbb{N}}
\def\Rit{\mathbb{R}}

\def\prob{\mathbb{P}}
\def\Prob{\mathbb{P}}

\begin{document}

\parindent=0.3in

\parindent=0.3in

\thispagestyle{empty}


\begin{center}
{\large\sc Uniform bias study and Bahadur representation \\
for local polynomial estimators of the conditional quantile function}\footnote{
This paper was started and completed when both authors were at Laboratoire de Statistique Th\'eorique et Appliqu\'ee, Universit\'e Pierre et Marie Curie, which support is gratefully acknowledged.
Financial Support from
the Department of Economics, Queen Mary University of London, is also gratefully
acknowledged. The authors would like to thank the participants of the Queen Mary Econometrics Reading Group, of the Berlin Quantile Regression Workshop, of the LSE Econometrics and Statistics Workshop as well as the Associate Editor and two anonymous referees whose careful readings, suggestions and comments have helped to improve the paper. All remaining errors are under our responsibility. This version corrects an error in the proof of Lemma B.2 which was pointed out by Zhongjun Qu but does not change the results of the published version.}



Emmanuel Guerre\footnote{
School of Economics and Finance, Queen Mary University of London, Great-Britain. \texttt{%
e.guerre@qmul.ac.uk }}\\[0pt]

\medskip

Camille Sabbah\footnote{Laboratoire EQUIPPE, Universit\'{e} Lille 3, France. {\tt  camille.sabbah@univ-lille3.fr}}
\\[0pt]
\vspace{0.5cm}


September 2014

\vspace{1cm}

\textbf{Abstract}\\[0pt]
\begin{minipage}[t]{12cm}
\begin{small}
This paper investigates the bias and the weak Bahadur representation of a local polynomial estimator of the conditional quantile function and its derivatives. The bias and Bahadur remainder term are studied uniformly with respect to the quantile level, the covariates and the smoothing parameter. The order of the local polynomial estimator can be higher than the differentiability order of the conditional quantile function. Applications of the results deal with global optimal consistency rates of the local polynomial quantile estimator, performance of random bandwidths and estimation of the conditional quantile density function. The latter allows to obtain a simple estimator of the conditional quantile function of the private values in a first price sealed bids auctions under the independent private values paradigm and risk neutrality. 

\vspace{1cm}


\vspace{.5cm}


\vspace{.5 cm}

\end{small}
\end{minipage}
\end{center}

\renewcommand{\thefootnote}{\arabic{footnote}} \setcounter{footnote}{0} %
\setcounter{page}{0} \newpage 

\newpage

\section{Introduction}
The conditional quantile function is a powerful tool to represent the dependence between two variables. Let $Q(\alpha|x)$, $\alpha$ in $(0,1)$, be the conditional quantile function of a univariate dependent variable $Y$ given $X=x$, where $X$ is the $d$ dimensional covariate, $Q(\alpha|x)=\inf \{ y : \Prob (Y\leq y |X=x) \geq \alpha \}$. Under fairly general conditions, the  L\'evy-Smirnov-Rosenblatt transformation ensures that there is a random variable $A$ independent of $X$ and uniform over $[0,1]$ such that
\begin{equation}
Y=Q(A|X).
\label{LSR}
\end{equation}
 In other words, the knowledge of the conditional quantile function allows to compute the impact on $Y$ of a shock on $X$ for any given $A$. The conditional quantile function is also central in the identification of the impact of such shocks or of more general parameters in nonseparable models in microeconometrics, see Chesher (2003), Chernozhukov and Hansen (2005), Holderlein and Mammen (2007) and Imbens and Newey (2009) to mention just a few. See also Firpo, Fortin and Lemieux (2009) or Rothe (2010) for an unconditional point of view when evaluating  distributional policy effects. Conditional quantile approaches can also be useful in industrial organization due to the important role played by increasing functions and the equivariance property of quantile function which states that $\Psi(Q(\alpha|x))$ is the conditional quantile function of $\Psi(Y)$ given $X$ provided $\Psi$ is an increasing transformation. See Haile, Hong and Shum (2003), Marmer and Shneyerov (2008) and below for the case of auctions. Echenique and Komunjer (2009) show the usefulness of a conditional quantile approach when analyzing general multiple equilibria economic models.

However, inference with the quantile representation (\ref{LSR}) is potentially difficult due to nonseparability. In a regression model $Y=m(X)+\varepsilon$ where $X$ and $\varepsilon$ are independent, the dependence between $Y$ and $X$ is summarized through the regression function $m(\cdot)$ and does not involve the unobserved noise $\varepsilon$. This contrasts with (\ref{LSR}) where the random variable $A$ may potentially change the shape of $x \mapsto Q(A|x)$. Hence, inference in (\ref{LSR}) should not focus on a particular value of the quantile level $\alpha$ but should consider instead all $\alpha$ in an interval $[\underline{\alpha},\overline{\alpha}]$ close enough to $[0,1]$, as recommended for instance in the case of the more constrained quantile regression model analyzed in Koenker (2005). In practice, this often leads to consider graphical representations of the estimated  curves $x \mapsto \widehat{Q} (\alpha|x)$ for various $\alpha$. A natural norm for evaluating  these estimated graphs is the uniform norm with respect to $\alpha$ and $x$, $\sup_{\alpha,x} \left| \widehat{Q} (\alpha|x)- Q(\alpha|x) \right|$.

The present paper contributes to this issue for local polynomial quantile estimators $\widehat{Q}_h (\alpha|x)$ which depends upon a bandwidth $h$. We study its bias uniformly in $\alpha$ and $x$ and derive a uniform Bahadur representation for $\widehat{Q}_h (\alpha|x)$ and its derivatives which holds in probability, that is a weak Bahadur representation. In few words, a Bahadur representation is an approximation of $\widehat{Q}_h (\alpha|x) - Q(\alpha|x)$ by a bias term plus a leading stochastic term up to remainder term with an explicit order. In our setup, uniformity is with respect to the level $\alpha$, the bandwidth $h$, and the covariate $x$, implying that our Bahadur representation is an important step for the study of $\sup_{\alpha,x} \left| \widehat{Q} (\alpha|x)- Q(\alpha|x) \right|$, see Proposition \ref{Randh} below. Various other interesting results also follow from our uniform results.

To be more specific, consider independent and identically observations $(X_1,Y_1),\ldots,(X_n,Y_n)$ with the same distribution than $(X,Y)$. Define, for $\alpha$ in $(0,1)$, the loss function
\begin{equation}
  \ell_{\alpha} ( q) = |q| + (2 \alpha -1) q
  = 2 q \left( \alpha- \ind (q \leq 0) \right)
  , \;
  q \mbox{ \rm in } \Rit,
\label{Lalpha}
\end{equation}
where $\Rit$ stands for the set of real numbers. 
It is well known that
\begin{equation}
Q \left(\alpha \left|x \right.\right)
=
\arg
\min_{q \in \Rit}
\esp
\left[
\ell_{\alpha} ( Y-q)
\left|X=x\right.
\right]
\label{Quant}
\end{equation}
is the conditional quantile of $Y$ given $X=x$. When $d=1$, the local polynomial estimator of order $p$ of $Q ( \alpha|x)$ is $\widehat{Q}_h (\alpha|x) = \widehat{ b}_0 (\alpha;h,x)$ where, for 
$
{\bf b}
=
\left(
b_0, \ldots, b_p
\right)^{T}
$,
\begin{equation}
\widehat{{\bf b}} (\alpha;h,x)
=
\arg
\min_{{\bf b} \in \Rit^{p+1}}
\sum_{i=1}^n
\ell_{\alpha}
\left(
Y_i
- 
b_0
-
b_1
\left(X_i-x\right)
-
\cdots
-
\frac{b_p}{p!}
\left(X_i-x\right)^p
\right)
K
\left(
\frac{X_i-x}{h}
\right)
.
\label{LPd1}
\end{equation}
In the expression above, $p!$ is the factorial $p \times (p-1) \times \cdots \times 1$, $K(\cdot)$ is a kernel function and $h$ is a smoothing parameter which goes to $0$ with the sample size. As detailed in Section 2 and studied throughout the paper, the local polynomial estimator $\widehat{Q}_h (\alpha|x)$ has a natural extension which covers the multivariate case $d>1$. As noted in Fan and Gijbels (1996, Chapter 5), the local polynomial estimator $\widehat{Q}_h (\alpha|x)$ is a modification of the Least Squares local polynomial estimator of a regression function which uses the square loss function in (\ref{LPd1}) instead of the loss function $\ell_{\alpha} (\cdot)$. A Taylor expansion
$$
Q (\alpha|X_i)
\simeq
Q(\alpha|x)
+
\frac{\partial Q (\alpha|x)}{\partial x}
\left(
X_i-x
\right)
+
\cdots
+
\frac{1}{p!}
\frac{\partial^p Q (\alpha|x)}{\partial x^p}
\left(
X_i-x
\right)^p
$$
suggests that $\widehat{ b}_1 (\alpha;h,x), \ldots ,\widehat{ b}_p (\alpha;h,x)$ estimate the partial derivatives 
$
\partial^r Q (\alpha|x) / \partial x^r
$, $r=1,\ldots,p$, 
provided $Q (\alpha |x)$ is smooth enough. 

Robust local polynomial estimation of a regression function and its derivatives, including quantile methods, has already been considered in many research articles. See in particular Tsybakov (1986) for optimal pointwise consistency rates, Fan (1992) for design adaptation, and Fan and Gijbels (1996) and Loader (1999) for a general overview. The present paper is perhaps more specifically related to Truong (1989), Chauduri (1991), Holderlein and Mammen (2009) and Kong, Linton and Xia (2010). Truong (1989) showed that local median estimators achieve the global optimal rates of Stone (1982) with respect to $L_{m}$ norms, $0<m\leq\infty$, for conditional quantile function satisfying a Lipschitz condition. Chauduri (1991) obtained a strong (that is which holds in an almost sure sense) Bahadur representation for the local polynomial quantile estimators when the kernel function $K(\cdot)$ of (\ref{LPd1}) is uniform. Hong (2003) extended this result to local polynomial robust M-estimation and more general kernels. The Bahadur representation of Chaudhuri (1991) is pointwise, that is holds for some prescribed $x$ and $\alpha$ and a given deterministic bandwidth $h\rightarrow 0$. As explained and illustrated in Kong et al. (2010), pointwise Bahadur representations are not sufficient for many applications including plug in estimation of conditional quantile functionals or marginal integration estimators. Hence Kong et al. (2010) derives a strong uniform Bahadur representation for robust local polynomial M-estimators for dependent observations. Here uniformity is with respect to the location variable $x$. For local polynomial quantile estimators of order $p=1$, Holderlein and Mammen (2009) considers uniformity with respect to $\alpha$ and $x$ but they just show that their remainder term is negligible in probability and does not obtain its order. 

In this work, we study the bias term and obtain the order in probability of the Bahadur remainder term uniformly in $\alpha$, $h$ and $x$ for local polynomial quantile estimators. A first contribution given in Theorem \ref{Bias} below deals with the study of the bias of local polynomial quantile estimators. Most of the literature has focused on the case where the order $p$ of the local polynomial is equal to the order of differentiability of $x \mapsto Q(\alpha|x)$, say $s$. This is somehow unrealistic since it amounts to assume that $s$ is known. Since the case where $p \leq s$ can be easily dealt with by ignoring higher order derivatives, we focus in the more interesting case where $p > s$, which has apparently not been considered in the statistical and econometric literature. As shown in  Corollary \ref{Optimal}, a local polynomial quantile estimator with $p > s$ still allows to estimate $Q(\alpha|x)$ with the optimal rate $n^{-s/(2s+d)}$ of Stone (1982). This suggests that local polynomial estimators using high order $p$ should be preferred since they allow to estimate in an optimal way a wider range of smooth conditional quantile functions. Another interesting conclusion of our bias study is that the additional local polynomial coefficients $\widehat{b}_v (\alpha;h,x)$, $v=s+1,\ldots,p$ can diverge and Proposition \ref{Sharp} describes a simple example where it indeed happens. Hence, in the local polynomial setup, a high value of $\widehat{b}_v (\alpha;h,x)$  may also correspond to a non smooth quantile function in which case a lower degree $p<v$ could have been used. 

Our uniform study of the Bahadur remainder term, namely Theorem \ref{Bahadur}, is the second main contribution of the paper. A third contribution builds on the fact that Theorems \ref{Bias} and \ref{Bahadur} hold uniformly with respect to $x$ in a compact inner subset of the support of $X$. Combining these results with a study of the stochastic part of the Bahadur representation allows us to show that the local polynomial quantile estimator achieves the global optimal rates of Stone (1982) for the $L_m$ and uniform norms provided the bandwidth goes to $0$ with an appropriate rate. This result, stated in Corollary \ref{Optimal}, is apparently new and extends Truong (1989) which is restricted to Lipshitz quantile functions, or Chauduri (1991) who considers pointwise optimality. A fourth contribution uses the fact that Theorems \ref{Bias} and \ref{Bahadur} hold uniformly with respect to $h$ in an interval $[\underline{h},\overline{h}]$. Proposition \ref{Randh} shows that a random bandwidth performs as well as its deterministic equivalent counterpart with respect to convergence rates of the uniform norm $\sup_{\alpha,x} \left|
\widehat{Q}_h (\alpha|x) - Q(\alpha|x) \right|$. Such a result gives a solid theoretical basis to Li and Racine (2008) suggestion of choosing the local polynomial bandwidth $h$ via a simpler cross validation procedure for the conditional cumulative distribution function. As mentioned earlier, uniformity with respect to $\alpha$ and $x$ is also useful for graphical representations of (\ref{LSR}).

A fifth contribution also exploits uniformity with respect to the quantile order $\alpha$.  Proposition~\ref{Quantdens} considers estimation of the conditional quantile density function
\begin{equation}
q (\alpha|x)
=
\frac{\partial Q(\alpha|x)}{\partial \alpha}=
\frac{1}{f \left( Q(\alpha|x)|x\right)}
.
\label{CQDF}
\end{equation}
As argued in Parzen (1979), the quantile density function $q (\alpha|x)$ or its inverse $1/q(\alpha|x)$ is a renormalization of the density function $f(y|x)$ which is well suited for statistical explanatory analysis. The function $q(\alpha|x)$ is also crucial for quantile based statistical inference. Indeed, the asymptotic variance of $\widehat{Q}_h (\alpha|x)$ is proportional to
$$
\frac{1}{n h} \frac{\alpha (1-\alpha)}{q^2 (\alpha|x) f(x)}
$$
where $f(\cdot)$ is the marginal density of $X$, see Fan and Gijbels (1996, p. 202). Hence estimating $q(\alpha|x)$ is useful to estimate the variance of $\widehat{Q}_h (\alpha|x)$. As noted in Guerre, Perrigne and Vuong (2009), the conditional quantile density function plays an important role in the identification of first-price sealed bids auction models. Under the independent private values paradigm and risk neutrality, the conditional quantile function of the private values $Q^v (\alpha|x)$ satisfies
$$
Q^v (\alpha|x) = Q^b (\alpha|x) +  \frac{\alpha q^b (\alpha|x)}{I-1} 
,
$$
where $Q^b (\alpha|x)$ and $q^b (\alpha|x)$ are the conditional quantile function and quantile density function of the bids. Hence estimating $Q^b (\alpha|x)$ and $q^b (\alpha|x)$ gives a straightforward estimation of the conditional quantile function of the private values $Q^v (\alpha|x)$ which is an alternative to the two steps approach of Guerre, Perrigne and Vuong (2000). See Haile et al. (2003) or Marmer and Shneyerov (2008) for a related estimation strategy.

There is however just a few references that address the estimation of $q(\alpha|x)$. For the related function $ q (\alpha|x) \partial F (Q(\alpha|x)|x) \partial x$, Lee and Lee (2008) uses a composition approach which nonparametrically estimates $\partial F(y|x) / \partial x$, $f(y|x)$ and $Q(\alpha|x) = F^{-1} (\alpha|x)$. Haile et al. (2003) and Marmer and Shneyerov (2008) proceeds similarly. Xiang (1995) proposes the estimator
$$
\frac{1}{h_q}\int \widehat{F}^{-1}  \left( \alpha + h_q a |x \right) d K_q(a),
$$
where $\widehat{F} (y|x)$ is a kernel estimator of the conditional cumulative distribution function, $K_q (\cdot)$ a probability distribution and $h_q$ a smoothing parameter. 
As argued in Fan and Gijbels (1996), local polynomial estimators may have better design adaptation properties than kernel ones.
Hence we propose to use the local polynomial $\widehat{Q}_h (\alpha|x)$ instead of the kernel $\widehat{F}^{-1} (\alpha|x)$. Thanks to uniformity with respect to $\alpha$ in Theorems \ref{Bias} and \ref{Bahadur}, the resulting conditional quantile density function estimator $\widehat{q} (\alpha|x)$ has a simple Bahadur representation which facilitates the study of its consistency rate, see Proposition \ref{Quantdens}.

The rest of the paper is organized as follows. The next section groups our main assumptions and notations and explained in particular how to extend (\ref{LPd1}) to multivariate covariates. Section 3 exposes our main results and Section 4 concludes the paper. The proofs of our statements are gathered in two appendices.

\setcounter{equation}{0} 
\section{Main assumptions and notations}

The definition (\ref{LPd1}) of $\widehat{Q}_h (\alpha|x)$ assumes that the covariate $X$ is univariate. In the multivariate case, we use a multivariate kernel function $K(z)=K(z_1,\ldots,z_d)$ but we restrict to an univariate bandwidth for the sake of simplicity. 
The univariate polynomial expansion  
$
b_0
+
b_1
\left(X_i-x\right)
+
\cdots
+
b_p
\left(X_i-x\right)^p/p!
$ is replaced by a multivariate counterpart as defined now. Let $\Nit$ be the set of natural integer numbers. For ${\bf v} = (v_1, \ldots, v_d)$ let 
$
|{\bf v}|=v_1+\cdots+v_d    
$
and let $P$ be the number of ${\bf v}$'s with $|{\bf v}| \leq p$. Then a generic expression for multivariate polynomial function of order $p$ is, for ${\bf b}$ in $\Rit^P$,
$$
{\bf U} (z)^{T} {\bf b}
=
\Sum_{{\bf v}; |{\bf v}| \leq p} b_{{\bf v}}\frac{z^{{\bf v}}}{{\bf v}!} ,
\mbox{ \rm where }
z^{{\bf v}} = z_1^{v_1} \times \cdots \times z_d^{v_d},
\; 
{\bf U} (z)^{T} 
=
\left(\frac{z^{{\bf v}}}{{\bf v}!} , |{\bf v}| \leq p \right),
$$
and ${\bf v}! = \Pi_{i=1}^d v_i !$.
In the expression above, the vectors ${\bf v}$ of $\Nit^d$ are ordered according to the lexicographic order. The multivariate version of the local polynomial estimator (\ref{LPd1}) is
\begin{eqnarray}
\widehat{{\bf b}} (\alpha;h,x)
& = &
\arg
\min_{{\bf b} \in \Rit^{P}}
\mathcal{L}_n
\left(
{\bf b}
;
\alpha,h,x
\right)
\mbox{ \rm with }
\label{LP}\\
\mathcal{L}_n
\left(
{\bf b}
;
\alpha,h,x
\right)
& = &
\frac{1}{nh^d}
\sum_{i=1}^n
\ell_{\alpha}
\left(
Y_i
- 
{\bf U} \left(X_i-x\right)^{T} {\bf b}
\right)
K
\left(
\frac{X_i-x}{h}
\right)
.
\nonumber
\end{eqnarray}
As in the univariate case, the entry $\widehat{b}_{{\bf 0}} (\alpha;h,x)=\widehat{Q}_h (\alpha|x)$ of $\widehat{{\bf b}} (\alpha;h,x)$ is an estimator of $Q(\alpha|x)$. The entry $\widehat{b}_{{\bf v}} (\alpha;h,x)$ can be viewed as an estimator of the partial derivative
$$
b_{{\bf v}} \left(\alpha|x\right)
=
\frac{
\partial^{|{\bf v}|} Q ( \alpha|x)
}{
\partial x_1^{v_1} \times \cdots \times \partial x_d^{v_d}
}
$$
provided this partial derivative exists. We shall consider later on the following 
Hölder class. Consider a subset $[\underline{\alpha},\overline{\alpha}]$ of $(0,1)$ over which $Q(\alpha|x)$ or its partial derivatives will be estimated. 
Let $\left\lfloor s \right\rfloor$ be the lowest integer part of 
$s$, i.e. $\left\lfloor s \right\rfloor$ is the unique integer number with
$
\left\lfloor
s
\right\rfloor
<
s
\leq
\left\lfloor
s
\right\rfloor
+
1
$.
Then $Q(\cdot|\cdot)$ is in $\mathcal{C} (L,s)$, $L,s>0$, if
\begin{enumerate}
\item 
for all $\alpha$ in $[\underline{\alpha},\overline{\alpha}]$, $x \mapsto Q(\alpha|x)$ is $\left\lfloor s \right\rfloor$-th continuously differentiable over the support $\mathcal{X}$ of $X$;
\item
for all ${\bf v}$ in $\Nit^d$ with $|{\bf v}| = \left\lfloor s \right\rfloor$,
    all $\alpha$ in $[\underline{\alpha},\overline{\alpha}]$, all $x$, $x'$ in $\mathcal{X}$,
\[
    \left|
      b_{{\bf v}} \left(\alpha|x\right)
      - 
      b_{{\bf v}} \left(\alpha|x'\right)
    \right|
    \leq
    L
    \left\|
      x - x'
    \right\|^{s - \lfloor s \rfloor }
\]
where $\| \cdot \|$ stands for the Euclidean norm.
\end{enumerate} 
Since the estimators 
$
\widehat{b}_{{\bf v}} \left(\alpha;h,x\right)
$
of the partial derivatives 
$b_{{\bf v}} \left(\alpha|x\right)$
converge with different rates, we use the diagonal standardization matrix
\[
{\bf H}
=
{\bf H}(h)
=
\Diag
\left(
h^{|{\bf v}|}
, {\bf v} \in \Nit^d, |{\bf v}| \leq p
\right)
.
\] 

It is well known that local polynomial estimation techniques apply at the boundaries. However we will focus on those $x$ which are in an inner subset $\mathcal{X}_0$ of the support $\mathcal{X}$ of $X$ to avoid technicalities. 
Our main assumptions are as follows. Let $\mathcal{B}\left(0,1\right)$ be the closed unit ball $\left\{z \in \Rit^d: \left\|z\right\| \leq 1 \right\}$.

\setcounter{hp}{23}
\begin{hp}
The distribution of $X$ has a probability density function $f (\cdot)$ with
respect to the Lebesgue measure, which is strictly positive and continuously differentiable over
the compact support $\mathcal{X}$ of $X$. The set $\mathcal{X}_0$ is a compact subset of the
interior of $\mathcal{X}$.
\label{X}
\end{hp}
\setcounter{hp}{5}
\begin{hp}
The cumulative distribution function $F(\cdot|\cdot)$ of $Y$ given $X$ has a continuous probability density function
$f (y|x)$ with respect to the Lebesgue measure, which is strictly positive for $y$ in $\mathbb{R}$ and $x$ in $\mathcal{X}$.
The partial derivative $\partial F(y|x) / \partial x$ is continuous  over $\Rit \times \mathcal{X}$.
There is a
$L_0 > 0$,
such that
\[
\left|
  f(y|x)
  -
  f(y'|x')
\right|
\leq
L_0
\left\|
  (x,y) - (x',y')
\right\|
\mbox{ \it for all $(x,y)$, $(x',y')$ of $\mathcal{X} \times \Rit$.}
\]
\label{F}
\end{hp}
\setcounter{hp}{10}
\begin{hp}
The nonnegative kernel function $K(\cdot)$ is Lipschitz over $\Rit^d$, has a compact support $\mathcal{K}$ and satisfies $\int K (z) dz = 1$. 
For some $\underline{K} > 0$, $K(z) \geq \underline{K}\; \ind\left( z \in \mathcal{B}\left(0,1\right)\right)$. 
The bandwidth is in $[\underline{h}_n,\overline{h}_n]$ with
$0 < \underline{h}_n \leq \overline{h}_n < \infty$, $\lim_{n \rightarrow \infty} \overline{h}_n = 0$ and  $\lim_{n \rightarrow \infty} (\log n) / (n \underline{h}^d_n ) =0$. 
\label{K}
\end{hp}
\noindent
Assumption \ref{X} is standard. Assumption \ref{F} ensures uniqueness of the conditional quantile $Q(\alpha|x)=F^{-1} (\alpha|x)$ in (\ref{Quant})
and existence of the quantile density function (\ref{CQDF}).
Assumption \ref{K} allows for a wide range of smoothing parameters $h \rightarrow 0$ in $[\underline{h}_n,\overline{h}_n]$. In the univariate case $d=1$, Hong (2003) restricts to bandwidths $h=O(n^{-1/(2p+3)})$, a condition which is not imposed here, and Chauduri assumes that $h$ has the exact order $n^{-1/(2p+d)}$. In the simpler context of univariate kernel regression, Einmahl and Mason (2005)  assumes   $h^d \geq C (\log n)/n$ to obtain uniform consistency so that Assumption \ref{K} is fairly general. 

\setcounter{equation}{0} 
\section{Bias study and Bahadur representation}
Applying standard parametric $M$-estimation theory as detailed in White (1994) or van der Vaart (1998) suggests that the local polynomial estimator $\widehat{{\bf b}} (\alpha;h,x)$ of (\ref{LP}) is an estimator of ${\bf b}^* (\alpha;h,x)$ with
\begin{equation}
{\bf b}^* (\alpha;h,x)
=
\arg
\min_{{\bf b} \in \Rit^{P}}
\esp
\left[
\ell_{\alpha}
\left(
Y
- 
{\bf U} \left(X-x\right)^{T} {\bf b}
\right)
K
\left(
\frac{X-x}{h}
\right)
\right]
.
\label{b*}
\end{equation}
In particular, $Q^*_h (\alpha|x) = b^*_{{\bf 0}} (\alpha;h,x)$ may differ from the true conditional quantile $Q(\alpha|x)$ due to a bias term
$Q^*_h (\alpha|x) - Q(\alpha|x)$. 
Studying this bias term can be done using the first-order condition 
$$
\frac{\partial}{\partial {\bf b}^T}
\esp
\left[
\ell_{\alpha}
\left(
Y
- 
{\bf U} \left(X-x\right)^{T} {\bf b}^* (\alpha;h,x)
\right)
K
\left(
\frac{X-x}{h}
\right)
\right]
=
0,
$$
and the Implicit Functions Theorem. This approach gives in particular the order of the difference between $b^*_{{\bf v}} (\alpha;h,x)$ and the ${\bf v}$th partial derivative $b_{{\bf v}} \left(\alpha|x\right)$ of $Q(\alpha|x)$ provided the partial derivative exists.
\begin{thm}
Assume that $Q(\cdot|\cdot)$ is in a H\"older class $\mathcal{C}(L,s)$ with
$\left\lfloor s \right\rfloor \leq p$. Then under Assumptions \ref{F}, \ref{K} and \ref{X} and provided $\overline{h}$ is small enough, there is a constant $C$ such that for all $|{\bf v}| \leq \lfloor s \rfloor$ and $n$ large enough,
\[
\sup_{
(\alpha,h,x)
   \in
      [\underline{\alpha},\overline{\alpha}]
    \times
      [\underline{h},\overline{h}]
     \times
       \mathcal{X}_0
     }
\left|
  \frac{
    b^{*}_{{\bf v}} (\alpha;h,x)
      -
    b_{{\bf v}} \left(\alpha|x\right)
       }
       {
       h^{s-|{\bf v}|}
       }
\right|
\leq C L
 .
\]
\label{Bias}
\end{thm}
\noindent
It follows that 
$
Q^*(\alpha|x) - Q(\alpha|x)
=
O(h^{s})
$ and more generally that
$$
b^{*}_{{\bf v}} (\alpha;h,x)
-
b_{{\bf v}} \left(\alpha|x\right)
=
O\left(h^{s-|{\bf v}|}\right)
$$
uniformly provided $|{\bf v}| \leq \lfloor s \rfloor$.
Since
$\left\lfloor s \right\rfloor \leq p$,
the bias order $h^{s-|{\bf v}|}$ is not affected by the order $p$ of the local polynomial estimator. This bias order is better than the bias order $h^{p-|{\bf v}|}$, $|{\bf v}| \leq p$, that would be achieved by suboptimal local polynomial estimators of lower order $p < \lfloor s \rfloor$.

The proof of Theorem \ref{Bias} establishes a slightly stronger result since it also gives the order of the coefficients $b^{*}_{{\bf v}} (\alpha;h,x)$ with $|{\bf v}| > \lfloor s \rfloor$ which correspond to partial derivatives that may not exist. Indeed, equation (\ref{Biascomplete}) of the proof of Theorem \ref{Bias} implies that
\begin{equation}
b^{*}_{{\bf v}} (\alpha;h,x)
=
O\left(h^{s-|{\bf v}|}\right)
\quad
\mbox{ \rm for $|{\bf v}|\geq s$}
\label{B*highv}
\end{equation}
uniformly in $(\alpha,h,x)
   \in
      [\underline{\alpha},\overline{\alpha}]
    \times
      [\underline{h},\overline{h}]
     \times
       \mathcal{X}_0$.
See also Loader (1999, Theorem 4.2) which gives a less precise 
$
b^{*}_{{\bf v}} (\alpha;h,x)
=
o\left(h^{-|{\bf v}|}\right)
$.
Hence the higher order polynomial coefficients $b^{*}_{{\bf v}} (\alpha;h,x)$, $|{\bf v}| > s$, may diverge when $h>0$. That this may be indeed the case can be seen on a simple regression example. Consider
\begin{equation}
Y= m(X) + \varepsilon,
\quad
m(x)
=
\left\{
\begin{array}{rl}
|x|^{1/2} & \mbox{ \rm if } x\geq0
\\
-|x|^{1/2} & \mbox{ \rm if } x<0
\end{array}
\right.
,
\label{Reg}
\end{equation}
where the $\mathcal{U} \left([-1,1]\right)$ random variable $X$  and the $\mathcal{N} \left(0,1\right)$ $\varepsilon$ are independent. Let $\Phi (\cdot)$ be the cumulative distribution function of the standard normal $\mathcal{N} \left(0,1\right)$. In this example, $Q(\alpha|x) = \Phi^{-1} (\alpha) + m(x)$ inherits of the smoothness properties of the regression function $m(\cdot)$. Note that the differential of $m(\cdot)$ at $x=0$ is infinite. It also follows that $Q (\alpha|x)$ 
is at best in an Hölder class $\mathcal{C}(L,1/2)$ since, for $L$ large enough,
$$
|m(x) -m(x')| \leq L \left|x-x'\right|^{1/2}
\quad
\mbox{ \rm for all $(x,x') \in [-1,1]^2$,}
$$
an inequality that cannot be improved by increasing the exponent $1/2$
as seen by taking $x=0$ and $x'\rightarrow 0$.
The next Proposition uses the behavior of $m(\cdot)$ at $x=0$ to show that the rate given in (\ref{B*highv}) is sharp.
\begin{prop}
Suppose that $(X,Y)$ satisfies (\ref{Reg}). Let ${\bf b}^*(\alpha;h,x)=\left(b_0^*(\alpha;h,x), b_1^*(\alpha;h,x)\right)^{T}$ from (\ref{b*}) be given by a local polynomial procedure of order 1. Then under Assumption \ref{K} and $\int z K(z) dz  = 0$,
$b_0^*(0.5;h,0)=m(0)+O(h^{1/2})$ and 
$b_1^*(0.5;h,0)$ diverges with the exact rate $h^{-1/2}$, 
$$
\lim_{h \rightarrow 0} h^{1/2}b_1^*(0.5;h,0) 
=
\frac{\int |z|^{3/2} K(z) d z}{\int z^{2} K(z) d z} \neq 0.
$$
\label{Sharp}
\end{prop}
\noindent
The divergence of $b_1^*(0.5;h,0)$ implies that the estimator $\widehat{b}_1(0.5;h,0)$ will diverge in probability. This recalls that observing a large $\widehat{b}_1(0.5;h,0)$ is not an argument for claiming that a local polynomial estimator of order $p=1$ should be used.

We now consider the stochastic terms $\widehat{Q}_h (\alpha|x) - Q^*_h(\alpha|x)$ and the rescaled 
$$
{\bf H} \left(\widehat{{\bf b}} (\alpha;h,x) - {\bf b}^* (\alpha;h,x)\right).
$$
Let us first introduce some additional notations. Local polynomial estimation builds on a order $p$ Taylor expansion of $Q(\alpha|x')$ with $x'$ in the vicinity of $x$. This Taylor expansion can be written as $Q(\alpha|x') \simeq {\bf U} (x'-x)^T {\bf b}_p (\alpha|x)$ where  ${\bf b}_p (\alpha|x)$ groups the partial derivatives of $Q(\alpha|x)$ with respect to $x$. Consider the following counterpart of
the Taylor approximation,  
\begin{equation}
Q^*( x';\alpha,h,x) = {\bf U} (x'-x)^T {\bf b}^* (\alpha,h,x)
\label{Q*}
\end{equation}
Define also ${\bf S}_i (\alpha;h,x)
=
{\bf S} (X_i,Y_i;\alpha,h,x)
$ and ${\bf J}_i ( \alpha;h,x)
=
{\bf J} (X_i; \alpha,h,x)$ with
\begin{equation}
{\bf S}_i (\alpha;h,x)
= 
2
\left\{
\ind
\left(
Y_i
\leq
Q^{*} (X_i;\alpha,h,x)
\right)
-
\alpha
\right\}
{\bf U}\left(\frac{X_i-x}{h}\right)
K
\left(\frac{X_i-x}{h}\right)
,
\label{Score}
\end{equation}
\begin{equation}
{\bf J}_i ( \alpha;h,x)
= 
2
f \left( Q^{*} (X_i;\alpha,h,x) \left| X_i \right. \right)
{\bf U} \left(\frac{X_i - x}{h}\right)
{\bf U} \left(\frac{X_i -x}{h}\right)^T
K
\left(\frac{X_i-x}{h}\right)
.
\label{J}
\end{equation}
Since
$$
{\bf U} \left(X_i - x\right)
= 
{\bf H}
{\bf U} \left(\frac{X_i - x}{h}\right) 
$$
and (\ref{Lalpha}) gives
\begin{eqnarray*}
\lefteqn{
\frac{\partial \ell_{\alpha}}{\partial {\bf b}^T}
\left(
Y_i
- 
{\bf U} \left(X_i-x\right)^{T} {\bf b}
\right)
K
\left(
\frac{X_i-x}{h}
\right)
}
&&
\\
& = & 
2
\left\{
\ind
\left(
Y_i \leq {\bf U} \left(X_i-x\right)^{T} {\bf b}
\right)
-
\alpha
\right\}
{\bf U} \left(X_i-x\right)
K
\left(
\frac{X_i-x}{h}
\right)
\end{eqnarray*}
almost everywhere,
the variables ${\bf S}_i (\alpha;h,x)$ satisfy
$$
\frac{\partial \mathcal{L}_n}{\partial {\bf b}^T}
\left({\bf b}^*(\alpha,h,x);\alpha,h,x\right)
=
\frac{{\bf H}}{n h^{d}}
\sum_{i=1}^n
{\bf S}_i (\alpha;h,x)
$$
almost everywhere. 
Hence $\sum_{i=1}^n
{\bf S}_i (\alpha;h,x)$ can be viewed as a score function term whereas $\sum_{i=1}^n
{\bf J}_i (\alpha;h,x)$ is actually similar to a second derivative of the objective function $\mathcal{L}_n$ although it is not twice differentiable.
Indeed, it can be shown that it admits a quadratic approximation  
with second-order derivatives
$$
{\bf H}
\left(
\frac{1}{n h^{d}}
\sum_{i=1}^n
{\bf J}_i (\alpha;h,x)
\right)
{\bf H}
.
$$
Classical results of White (1994) or van der Vaart (1998) for parametric estimation suggests that a candidate approximation for  
$\widehat{{\bf b}} (\alpha;h,x) - {\bf b}^* (\alpha;h,x)$ is
\begin{eqnarray*}
\lefteqn{
-\left( 
{\bf H}
\left(
\frac{1}{n h^{d}}
\sum_{i=1}^n
{\bf J}_i (\alpha;h,x)
\right)
{\bf H}
\right)^{-1}
\frac{{\bf H}}{n h^{d}}
\sum_{i=1}^n
{\bf S}_i (\alpha;h,x)
}
&&
\\
&=&
-
{\bf H}^{-1}
\left( 
\frac{1}{n h^{d}}
\sum_{i=1}^n
{\bf J}_i (\alpha;h,x)
\right)^{-1}
\frac{{1}}{n h^{d}}
\sum_{i=1}^n
{\bf S}_i (\alpha;h,x)
 .
\end{eqnarray*}
Hence the rescaled 
$
\left(n h^d\right)^{1/2} {\bf H}
\left(\widehat{{\bf b}} (\alpha;h,x) - {\bf b}^* (\alpha;h,x)\right)
$
is expected to be close to
\begin{equation}
\beta_n (\alpha;h,x)
=
-
\left(
  \frac{1}{nh^d}
  \Sum_{i=1}^{n}
    {\bf J}_i(\alpha;h,x)
\right)^{-1}
\frac{1}{\left(nh^d\right)^{1/2}}
  \Sum_{i=1}^{n} {\bf S}_i(\alpha;h,x).
\label{Betan}
\end{equation}
$\Sum_{i=1}^{n}{\bf J}_i(\alpha;h,x)/(nh^d)$ is similar to a Kernel regression estimator and obeys a Law of Large Numbers for triangular array which ensures that this matrix is asymptotically close to
$$
2
f \left( Q^{*} (x;\alpha,h,x) \left| x \right. \right)
\int
{\bf U} \left(t\right)
{\bf U} \left(t\right)^T
K
\left(t\right)
dt.
$$
Since this matrix is symmetric positive definite, the inverse in (\ref{Betan}) exists with a probability tending to 1. The term $\Sum_{i=1}^{n}{\bf S}_i(\alpha;h,x)/(nh^d)^{1/2}$ has a similar kernel structure but with centered 
${\bf S}_i(\alpha;h,x)$, see (\ref{FOC}) in Lemma \ref{Biais1} of Appendix A. Hence   
$\Sum_{i=1}^{n}{\bf S}_i(\alpha;h,x)/(nh^d)^{1/2}$ 
satisfies a pointwise Central Limit Theorem, as $\beta_n (\alpha;h,x)$. Hence 
$
\left(n h^d\right)^{1/2} {\bf H}
\left(\widehat{{\bf b}} (\alpha;h,x) - {\bf b}^* (\alpha;h,x)\right)
$
should also be asymptotically Gaussian provided the so called Bahadur error term
\begin{equation}
{\bf E}_n (\alpha;h,x)
=
\left(n h^d\right)^{1/2} {\bf H}
\left(\widehat{{\bf b}} (\alpha;h,x) - {\bf b}^* (\alpha;h,x)\right)
-
\beta_n (\alpha;h,x).
\label{En}
\end{equation}
is asymptotically negligible pointwisely. But transposing the various uniform results established in the Appendices for the leading term $\beta_n (\alpha;h,x)$ of the expansion of  $\left(n h^d\right)^{1/2} {\bf H}
\left(\widehat{{\bf b}} (\alpha;h,x) - {\bf b}^* (\alpha;h,x)\right)$ requests a uniform study of ${\bf E}_n (\alpha;h,x)$.

Techniques to study ${\bf E}_n (\alpha;h,x)$ for a fixed argument $\alpha$, $h$ and $x$ are given in Hjort and Pollard (1993). See also Fan, Heckman and Wand (1995, p.143) or Fan and Gijbels (1996, p.210). In our uniform setup, obtaining an uniform order for ${\bf E}_n (\alpha;h,x)$ is performed using a preliminary uniform study of a stochastic process we introduce now. Define first
\begin{eqnarray*}
\lefteqn{
\mathbb{L}_{1n} (\beta;\alpha,h,x)
}
&&
\\
& = &
nh^d
\left\{
\mathcal{L}_n
\left(
{\bf b}^* (\alpha;h,x)
+
\frac{H^{-1}\beta}{\left(n h^d\right)^{1/2}} 
;
\alpha,h,x
\right)
-
\mathcal{L}_n
\left(
{\bf b}^* (\alpha;h,x)
;
\alpha,h,x
\right)
\right\}
\\
& = &
\sum_{i=1}^n
\left\{
\ell_{\alpha}
\left(
  Y_i - Q^* (X_i ; \alpha,h,x)
  -
  \frac{{\bf U} \left(\frac{X_i-x}{h}\right)^T}{\left(n h^d\right)^{1/2}}
  \beta
\right)
-
\ell_{\alpha}
\left(
  Y_i - Q^* (X_i ; \alpha,h,x)
\right)
\right\}
K
\left(\frac{X_i-x}{h}\right),
\end{eqnarray*}
which is such that
$$
\left(n h^d\right)^{1/2} {\bf H}
\left(\widehat{{\bf b}} (\alpha;h,x) - {\bf b}^* (\alpha;h,x)\right)
=
\arg
\min_{\beta}
\mathbb{L}_{1n} (\beta;\alpha,h,x)
.
$$
It then follows from (\ref{En}) that
\begin{eqnarray}
{\bf E}_n (\alpha;h,x)
& = &
\arg \min_{\epsilon}
\mathbb{L}_{n}
\left(\beta_n(\alpha;h,x),\epsilon; \alpha;h,x\right)
\mbox{ \rm where }
\nonumber \\
\mathbb{L}_{n}
\left(\beta,\epsilon; \alpha;h,x\right)
& = &
\mathbb{L}_{1n} (\beta+\epsilon;\alpha,h,x)
-
\mathbb{L}_{1n} (\beta;\alpha,h,x)
.
\label{LL}
\end{eqnarray}
Hence the stochastic process $\mathbb{L}_n$ plays a central role in our analysis. Especially useful is the decomposition
$$
\mathbb{L}_{n}
\left(\beta,\epsilon; \alpha;h,x\right)
=
\mathbb{L}_{n}^{0}
\left(\beta,\epsilon; \alpha;h,x\right)
+
\mathbb{R}_{n}
\left(\beta,\epsilon; \alpha;h,x\right)
$$
where $\mathbb{L}_{n}^{0}$ is the quadratic approximation of
$\mathbb{L}_{n}$,
\begin{eqnarray}
\mathbb{L}_{n}^{0}
\left(\beta,\epsilon; \alpha;h,x\right)
& = &
\frac{1}{\left(n h^d\right)^{1/2}}
\sum_{i=1}^{n}
{\bf S}_i (\alpha;h,x)^{T} \left(\beta+\epsilon\right)
+
\frac{1}{2}
\left(\beta+\epsilon\right)^{T}
\left(
\frac{1}{n h^d}
\sum_{i=1}^n
{\bf J}_i (\alpha;h,x)
\right) 
\left(\beta+\epsilon\right)
\nonumber \\
&&
-
\frac{1}{\left(n h^d\right)^{1/2}}
\sum_{i=1}^{n}
{\bf S}_i (\alpha;h,x)^{T} \beta
+
\frac{1}{2}
\beta^{T}
\left(
\frac{1}{n h^d}
\sum_{i=1}^n
{\bf J}_i (\alpha;h,x)
\right) 
\beta
\nonumber
\\
& = &
\frac{1}{\left(n h^d\right)^{1/2}}
\sum_{i=1}^{n}
{\bf S}_i (\alpha;h,x)^{T} \epsilon
+
\frac{1}{2}
\epsilon^{T}
\left(
\frac{1}{n h^d}
\sum_{i=1}^n
{\bf J}_i (\alpha;h,x)
\right) 
\left(\epsilon + 2 \beta \right),
\label{LL0}
\end{eqnarray}
and $\mathbb{R}_{n}$ is a remainder term.
As in the expression above (\ref{LL}) for
${\bf E}_n (\alpha;h,x)$,
the variable $\beta$ above in (\ref{LL0}) will be taken equal to
$\beta_n(\alpha;h,x)$ in the proof of Theorem \ref{Bahadur} below.
As noted in the quadratic approximation lemma of Fan et al. (1995, p.148) in the pointwise case, the order of ${\bf E}_n (\alpha;h,x)$ is driven by the order of $\mathbb{R}_{n}$. The proof of the next Theorem relies on an uniform study of $\mathbb{R}_{n}$ based on a maximal inequality under bracketing entropy conditions from Massart (2007), see the proof of Proposition A.1. This maximal inequality plays here the role of the Bernstein inequality  used in the pointwise framework of Hong (2003). 
\begin{thm}
Under Assumptions \ref{F}, \ref{K} and \ref{X},
\[
\sup_{
(\alpha,h,x)
\in
[ \underline{\alpha},\overline{\alpha} ]
\times
[ \underline{h}, \overline{h}]
\times
\mathcal{X}_0
}
\left\|
{\bf E}_n (\alpha;h,x)
\right\|
=
O_{\prob}
\left(
\frac{
    \log^3
     \left(
       n
     \right)
     }
     {
     n \underline{h}^d
     }
\right)^{1/4}
.
\]
\label{Bahadur}
\end{thm}
\noindent
In the case where the lower and upper bandwidths $\underline{h}$ and $\overline{h}$ have the same order, Theorem \ref{Bahadur} gives uniformly in $h$ in $[\underline{h},\overline{h}]$, $\alpha$ and $x$,
$$
\widehat{Q}_h (\alpha|x)
= 
Q^*_h (\alpha|x)
+
\frac{{\bf e}_0^{T} \beta_n (\alpha;h,x)}{\left( nh^d \right)^{1/2}}
+
O_{\prob}
\left(
\frac{\log n}{nh^d}
\right)^{3/4},
$$
where ${\bf e}_0$ is the first vector of the canonical basis of $\Rit^P$, which  first coordinate is equal to $1$ and the other ones are equal to $0$. For $h$ of order $n^{-1/(2p+d)}$ as studied in Chauduri (1991, Theorem 3.2), the order of the remainder term is 
$n^{-3p/(2(2p+d))} \log^{3/4} n$ as found by this author. When $d=1$, Hong (2003) obtains the better order $(\log \log n /(nh))^{-3/4}$ but his Bahadur representation only holds pointwisely in $\alpha$ and $x$. It can be conjectured that the order $(\log n /(nh^d))^{-3/4}$ is optimal for Bahadur expansion holding uniformly with respect to $x$.

For higher order partial derivatives, Theorem \ref{Bahadur} yields
$$
\widehat{b}_{{\bf v}} (\alpha;h,x)
=
b^* (\alpha;h,x)
+
\frac{{\bf e}_{\bf v}^{T} \beta_n (\alpha;h,x)}{ \left( nh^d\right)^{1/2} h^{|{\bf v}|}}
+
\frac{
1}
{h^{|{\bf v}|}}
O_{\prob}
\left(
\frac{\log n}{nh^d}
\right)^{3/4}
,
$$ 
where the ${\bf v}$th entry of ${\bf e}_{\bf v}$ is $1$ and the other are $0$, see also Hong (2003) for a pointwise version of this expansion and Kong et al. (2010) for a version which is uniform with respect to $x$.  Such expansion can be used to study the pointwise asymptotic normality of the local polynomial quantile estimator. Combining this Bahadur representation with the bias study of Theorem \ref{Bias} gives a global rate result which is apparently new. The next Corollary extends the study of local medians in Truong (1989).
\begin{cor}
Assume that $Q(\alpha|x)$ is in $\mathcal{C} (L,s)$ for some $\left\lfloor s\right\rfloor
\leq p$. Suppose that Assumptions \ref{F}, \ref{K} and \ref{X} hold. Then for all partial derivative order ${\bf v}$ with $\left|{\bf v}\right|\leq \left\lfloor s\right\rfloor$ and all $\alpha$ in $[\underline{\alpha},\overline{\alpha}]$,
\begin{enumerate}
\item
$
\left(
\int_{\mathcal{X}_0}
\left|
\widehat{b}_{\bf v}
\left(\alpha;h,x)\right)
-
b_{\bf v} (\alpha|x)
\right|^m
d x
\right)^{1/m}
=
O_{\prob}
\left(
\frac{1}{n}
\right)^{\frac{s-|{\bf v}|}{2s+d}}
$
for any finite $m>0$
provided $h$ is asymptotically proportional to $n^{-\frac{1}{2s+d}}$;
\item
$
\sup_{x \in \mathcal{X}_0}
\left|
\widehat{b}_{\bf v}
\left(\alpha;h,x\right)
-
b_{\bf v} (\alpha|x)
\right|
=
O_{\prob}
\left(
\frac{\log n}{n}
\right)^{\frac{s-|{\bf v}|}{2s+d}}
$
if $h$ is asymptotically proportional to $\left(\frac{\log n}{n} \right)^{\frac{1}{2s+d}}$.
\end{enumerate}
\label{Optimal}
\end{cor}
\noindent
Since the $b_{\bf v} (\alpha|x)$ are estimators of the partial derivatives of $m(x)$ in a regression model as (\ref{Reg}), 
It follows from Stone (1982) that the global rates derived in Corollary \ref{Optimal} are optimal in a minimax sense.

A second application builds on the uniformity with respect to the bandwidth $h$ of our Bahadur representation. The next Proposition allows for data-driven bandwidths. Observe that it also deals with the uniform norm
$
\sup_{(\alpha,x) \in [\underline{\alpha},\overline{\alpha}] \times \mathcal{X}_0}
\left|
\widehat{Q}_{h} (\alpha|x)
-
Q (\alpha | x)
\right|
$
which evaluates the estimated curves $(\alpha,x) \mapsto \widehat{Q}_{h} (\alpha|x)$  used in empirical graphic illustrations of  (\ref{LSR}). 
\begin{prop}
Consider a random bandwidth $\widehat{h}_n$ such that $\widehat{h}_n = O_{\prob} (h_n)$ and $1/\widehat{h}_n = O_{\prob} (1/h_n)$ where $h_n$ is a deterministic sequence satisfying 
$h_n = o(1)$ and $\lim_{n\rightarrow\infty} (\log n) / (nh_n^d) =0$. Suppose that Assumption \ref{K}, \ref{F} and \ref{X} hold and that $Q(\alpha|x)$ is in $\mathcal{C} (L,s)$. Then for any 
${\bf v}$ with $\left|{\bf v}\right|\leq \left\lfloor s\right\rfloor$,
\begin{eqnarray*}
\sup_{(\alpha,x) \in [\underline{\alpha},\overline{\alpha}] \times \mathcal{X}_0}
\left|
\widehat{b}_{\bf v} (\alpha; \widehat{h}_n,x)
-
b_{\bf v} (\alpha | x)
\right|
& = &
h_n^{-\left|{\bf v}\right|}
O_{\prob}
\left(
h_n^{s}
+
\left(
\frac{\log n}{ nh_n^d }
\right)^{1/2}
\right)
.
\end{eqnarray*}
\label{Randh}
\end{prop}
\noindent
In particular if the exact order of $\widehat{h}_n$ is $(\log (n)/n)^{1/(2s+d)}$ in probability,
$
\sup_{x \in \mathcal{X}_0}
 \left|
\widehat{b}_{\bf v} (\alpha; \widehat{h},x)
-
b_{\bf v} (\alpha | x)
\right| 
$
has the optimal order 
$(\log (n)/n)^{(s-\left|{\bf v}\right|)/(2s+d)}$ 
of Corollary \ref{Optimal}-(ii). It is likely that an $L_m$ version of Proposition \ref{Randh} holds but it is slightly longer to prove. Proposition \ref{Randh} can be for instance fruitfully applied to cross-validated bandwidths for the conditional cumulative distribution as proposed by Li and Racine (2008).

Our last application builds on the fact that Theorems \ref{Bias} and \ref{Bahadur} hold uniformly with respect to the quantile order $\alpha$. This application concerns estimation of the conditional quantile density function (\ref{CQDF}).
The considered estimator of $q(\alpha|x)$ is a conditional version of the Parzen (1979) convolution estimator, 
\begin{equation}
\widehat{q} (\alpha|x)
=
\frac{1}{h_q}
\int
\widehat{Q}_{h} (a|x)
d
K_q
\left(\frac{a-\alpha}{h_q}\right)
=
\frac{1}{h_q}
\int
\widehat{Q}_{h} (\alpha+h_q t|x)
d
K_q
\left(t\right)
,
\label{Hatq}
\end{equation}
see also Xiang (1995).
In the expression above, $h_q>0$ is a bandwidth and $K_q (\cdot)$ is a signed measure over $\Rit$ such that
$$
\int
d
K_q
\left(t\right)
=
0,
\quad
\int
t
d
K_q
\left(t\right)
=
1.$$ 
In particular, if $K_q (\cdot)$ has a Lebesgue derivative $d K_q (t) =K_q' (t)dt$, substituting in (\ref{Hatq}) gives
$$
\widehat{q} (\alpha|x)
=
\frac{1}{h_q}
\int
\widehat{Q}_h (\alpha+h_qt|x)
K_q'\left(t\right)
d t
.
$$
Computing these integrals may request intensive numerical steps so that the resulting estimator may be difficult to implement in practice. A more realistic estimator  uses a discrete measure $K_q (\cdot)$ in (\ref{Hatq}). If $K_q (\cdot)$ is a linear combination of Dirac masses at $t_j$ with weights $\kappa_j$, $j=1,\ldots,J$, the resulting estimator
$$
\widehat{q} (\alpha|x)
=
\frac{1}{h_q}
\sum_{j=1}^J
\kappa_j
\widehat{Q}_h (\alpha+h_q t_j|x),
\quad
\sum_{j=1}^J
\kappa_j
=
0
\mbox{ \rm and }
\sum_{j=1}^J
t_j
\kappa_j
=
1
,
$$
may be indeed simpler to compute. Note that this includes the well known numerical derivatives
$$
\frac{
\widehat{Q}_h (\alpha+h_q |x)
-
\widehat{Q}_h (\alpha |x)
}{h_q}
,
\quad
\frac{
\widehat{Q}_h (\alpha |x)
-
\widehat{Q}_h (\alpha-h_q |x)
}{h_q}
\mbox{ \rm and }
\frac{
\widehat{Q}_h (\alpha+h_q |x)
-
\widehat{Q}_h (\alpha-h_q |x)
}{2h_q}
.
$$
To study the bias of $\widehat{q} (\alpha|x)$, we strengthen the definition of the smoothness class $\mathcal{C} (L,s)$ as follows. $Q (\alpha|x)$ is in $\mathcal{C}_q (L,s)$ if 
\begin{enumerate}
	\item 
	$Q (\alpha|x)$ is in $\mathcal{C} (L,s+1)$;
	\item
	For each $x$ in $\mathcal{X}$, $\alpha \in [\underline{\alpha},\overline{\alpha}] \mapsto q(\alpha|x)$ is $\left\lfloor s\right\rfloor$th differentiable;
	\item 
	For each $x$ in $\mathcal{X}$ and all $(\alpha,\alpha') \in [\underline{\alpha},\overline{\alpha}]^2$
	$$
	\left|
	\frac{\partial^{\left\lfloor s\right\rfloor} q(\alpha|x)}{\partial \alpha^{\left\lfloor s\right\rfloor}}
	-
	\frac{\partial^{\left\lfloor s\right\rfloor} q(\alpha'|x)}{\partial \alpha^{\left\lfloor s\right\rfloor}}
	\right|
	\leq
	L
	\left|
	\alpha - \alpha'
	\right|^{s-\left\lfloor s\right\rfloor}
	.
	$$
\end{enumerate}
We shall assume in addition that $K_q(\cdot)$ has a compact support and satisfies the additional conditions
$$
\int t^j d K_q (t) = 0,
\quad
j=1, \ldots, \left\lfloor s\right\rfloor,
\quad
\int  \left|d K_q (t)\right| < \infty.
$$
\begin{prop}
Assume that $Q(\alpha|x)$ is in $\mathcal{C}_q (L,s)$ and $\left\lfloor s+1\right\rfloor \leq p$.  Suppose that Assumptions \ref{K}, \ref{F} and \ref{X} hold with
$h=O(h_q)$, $h_q \rightarrow 0$ and $ (\log n)/(nh^d) \rightarrow 0$. Then for any $x$ in $\mathcal{X}_0$ and $\alpha$ in $(\underline{\alpha},\overline{\alpha})$,
$$
\widehat{q}(\alpha|x)
=
q(\alpha|x)
+
O_{\prob}
\left(
h_q^s
+
\frac{1}{\left(n h^d h_q\right)^{1/2}}
\right)
+
\frac{\log^{3/4} n}{\left( nh^d h_q^2\right)^{1/4}}
O_{\prob}
\left(
\frac{1}{\left( n h_q h^d \right)^{1/2}}
\right).
$$
\label{Quantdens}
\end{prop}
\noindent
Taking $h_q$ and $h$ of the same order is the optimal choice for the order of $h$ in the expansion of Proposition \ref{Quantdens}. This gives 
$$
\widehat{q}(\alpha|x)
=
q(\alpha|x)
+
O_{\prob}
\left(
h^s
+
\frac{1}{\left(n h^{d+1} \right)^{1/2}}
\right)
+
\frac{\log^{3/4} n}{\left( nh^{d+2} \right)^{1/4}}
O_{\prob}
\left(
\frac{1}{\left( n h^{d+1} \right)^{1/2}}
\right).
$$
The item $\left( \log^{3/4} n\right)\left( nh^{d+2} \right)^{-1/4} O_{\prob}
\left(
\left( n h^{d+1} \right)^{-1/2}
\right)
$
is given by the Bahadur error term $\mathbf{E}_n (\alpha;h,x)$ of Theorem \ref{Bahadur}. The other item, $O_{\prob} \left(h^{s}+(nh^{d+1})^{-1/2}\right)$, can be viewed as a bias variance decomposition component. The latter is the leading term of the expansion provided 
$n h^{d+2} \rightarrow \infty$, a condition also used in Lee and Lee (2008) when $d=1$. In this case, the optimal order for $h$ is $n^{-1/(2s+d+1)}$, which is such that  
$n h^{d+2} \rightarrow \infty$ provided $s>1/2$. In this case, the optimal rate for pointwise estimation of $q(\alpha|x)$ is $n^{-s/(2s+d+1)}$ which, as expected from (\ref{CQDF}), coincides with the optimal rate for pointwise estimation of $f(y|x)$. 

\section{Final remarks}

This paper has investigated the bias and the Bahadur representation of a local polynomial estimator of the conditional quantile function and its derivatives. Compared to the existing literature, a distinctive feature is that the bias and Bahadur remainder term are studied uniformly with respect to the quantile level, the covariates and the smoothing parameter, extending so Chauduri (1991) and Kong et al. (2010). Our framework also considers the case where the order of the local polynomial estimator $p$ is higher than the order of differentiability $s$ of the conditional quantile function. An interesting consequence of our bias study is that using a local polynomial estimator of order $p \geq s$  does not affect its rate optimality.  

Our uniform study of the bias and of the Bahadur remainder term are applied to derive the global rate optimality of the local polynomial estimators of the conditional quantile function and its derivatives with respect to $L_m$ norms, $0<m\leq \infty$ provided the bandwidth goes to $0$ with an appropriate rate. This extends Truong (1989) who states a similar result for local medians and under a rather strong Lipschitz condition for the conditional quantile function. Another application deals with the performance of randomly selected bandwidths that are shown to perform as well as their deterministic equivalent in term of consistency rates in uniform norm. Our framework is flexible enough to be adapted to other global norms. This new result is especially useful in view of Li and Racine (2008) suggestion of implementing local polynomial quantile estimation with a data-driven bandwidth given by a cross validation criterion for the conditional cumulative distribution function. A last application to nonparametric estimation of the quantile density function can be useful for confidence intervals and in Econometrics of Auctions where the conditional quantile density function plays an important role.

Our uniform results can also be useful for other studies. For instance an issue far beyond the scope of the present paper is the choice of the local polynomial order $p$. Local polynomial quantile estimation can be implemented using a large $p$, possibly growing with the sample size. This would allow to estimate very smooth conditional quantile function with a small bias although it may inflate the asymptotic variance of the resulting estimator. Another approach would be to use a data-driven local polynomial order $p$. Such a problem is very close to the issue of choosing the order of the kernel when estimating a regression or a probability density function. The latter can be addressed following the recent adaptive approach of Goldenshluger and Lespki (2008,2009) which gives a data-driven choice of the kernel and bandwidth  in the context of the continuous time white noise model. Our uniform Bahadur representation is a preliminary step that can be useful to extend their results to local polynomial quantile estimation.

{\small
\renewcommand{\baselinestretch}{1.2} 
\setcounter{lem}{0} 
\setcounter{prop}{0}
\setcounter{thm}{0} 
\renewcommand{\thelem}{A.\arabic{lem}} 
\renewcommand{\theprop}{A.\arabic{prop}} 
\setcounter{equation}{0} \setcounter{subsection}{0} 
\renewcommand{\theequation}{A.\arabic{equation}} 
\renewcommand{\thesubsection}{A.\arabic{subsection}}
\section*{Appendix A: Proofs of main results}

Appendix A groups the proofs of Theorems \ref{Bias} and \ref{Bahadur}, Propositions \ref{Sharp}, \ref{Randh} and \ref{Quantdens}, and Corollary \ref{Optimal}. The proofs of intermediary results used to prove these main results are grouped in Appendix B.

We first introduce some additional notations. Sequences $\{a_n\}$ and $\{b_n \}$  satisfy
$a_n \asymp b_n$ if $|a_n|/C \leq |b_n| \leq C|a_n|$ for some $C>0$ and $n$ large enough.
Recall that $\| \cdot \|$ is the Euclidean norm and $\mathcal{B} (0,1)=\{z;\|z\| \leq 1\}$. Let $\succ$ be the usual order for symmetric matrices, that is ${\bf A_1} \succ {\bf A_2}$ if and only if ${\bf A_1}- {\bf A_2}$ is a non-negative symmetric matrix.
If ${\bf A}$ is a symmetric matrix, $
\left\|{\bf A}\right\|
=
\sup_{{\bf u} \in \mathcal{B} (0,1)}
\left\|{\bf A}{\bf u}\right\|
= 
\sup_{{\bf u} \in \mathcal{B} (0,1)} |{\bf u}^T{\bf A} {\bf u}|
$ is the largest eigenvalue in absolute value of ${\bf A}$. This norm is such that $\| {\bf A}{\bf  B} \| \leq \|{\bf A} \| \|{\bf B}\|$ for any matrix or vector ${\bf B}$.
Denote by $\| \cdot \|_{\infty}$ the uniform norm, i.e.
$
\left\| f (\cdot|\cdot) \right\|_{\infty}
=
\sup_{(x,y) \in \Rit^d \times \Rit}
| f(y|x) |
$.
We use the abbreviation $\theta = (\alpha,h,x)$. In particular, $Q^{*}(x';\theta)$, ${\bf S}_i (\theta)$ and ${\bf J}_i (\theta)$ stand for $Q^{*}(x';\alpha,h,x)$, ${\bf S} (X_i,Y_i;\alpha,h,x)$ and ${\bf J}(X_i;\alpha,h,x)$, see equations (\ref{Q*}), (\ref{Score}) and (\ref{J}).
We abbreviate $\underline{h}_n$ and $\overline{h}_n$ into $\underline{h}$ and $\overline{h}$. Define
\[
\Theta^0
=
[\underline{\alpha},\overline{\alpha}]
\times
[0 ,\overline{h}]
\times
\mathcal{X}_0
\; , \;
\Theta^1
=
[\underline{\alpha},\overline{\alpha}]
\times
[\underline{h} ,\overline{h}]
\times
\mathcal{X}_0
 ,
\]
where $\mathcal{X}_0$ is as in Assumption \ref{X} and $[\underline{\alpha} , \overline{\alpha}] \subset (0,1)$ is as in the definition of the smoothness class $\mathcal{C} (L,s)$.
For $\mathcal{L}_n 
\left({\bf b}; \alpha,h,x\right)=\mathcal{L}_n 
\left({\bf b};\theta\right)$ as in (\ref{LP}), define 
$$
\mathcal{L}
\left({\bf b};\theta \right)
= 
\esp
\left[
\mathcal{L}_n 
\left({\bf b};\theta\right)
\right]
=
\frac{1}{h^d}
\esp
\left[
\left\{
\ell_{\alpha}
\left(
Y
- 
{\bf U} \left(X-x\right)^{T} {\bf b}
\right)
-
\ell_{\alpha} \left(Y\right)
\right\}
K
\left(
\frac{X-x}{h}
\right)
\right]
.
$$
We also use $K_h (z) = K(z/h)$.
It is convenient to change ${\bf b}$ into its standardization
${\bf B} = {\bf H b}$ and to define
$
\widehat{{\bf B}} (\theta) = {\bf H} \widehat{{\bf b}} (\theta)
$ and
$
{\bf B}^*(\theta) = {\bf H} {\bf b}^* (\theta)$.
Absolute constants are denoted by the generic letter $C$ and may vary from line to line.

The following argument is used systemically. Recall that $\mathcal{X}_0$ is an inner subset of the compact $\mathcal{X}$ under Assumption \ref{X}. Hence for any $(x,h) \in \mathcal{X}_0 \times \mathcal{K}$, $x+h z$ is in $\mathcal{X}$ under Assumption \ref{K} provided $\overline{h}$ is small enough.

The next lemma is used in the proof of Theorems \ref{Bias} and \ref{Bahadur}. Its proof  is given in Appendix B with the proof of the other intermediary results.
\begin{lem}
Under Assumption \ref{F}, \ref{K} and \ref{X}, we have for $\overline{h}$ small enough,
\begin{enumerate}
\item
${\bf b}^{*}(\theta)$ exists and is unique for all $\theta$ in $\Theta^0$. 
\item
${\bf B}^{*}(\theta) = {\bf H}{\bf b}^{*}(\theta)$ satisfies 
\begin{eqnarray}
&&
\esp \left[{\bf S}_i \left(\theta\right)\right]
 = 
\int
\left\{
F
\left(
{\bf U}\left(z\right)^T {\bf B}^* (\theta)| x+hz 
\right)
-
F \left(Q(\alpha|x+hz)\left|x+hz\right.\right)
\right\}
   f(x+hz)
   {\bf U}(z)
   K(z)
 dz
=
 0,
\label{FOC}
\\
&&
\lim_{\overline{h} \rightarrow 0}
\sup_{\theta \in \Theta^0}
  \left\|
    {\bf B}^{*}(\theta)
      -
    {\bf B}^{*}(\alpha;0,x)
  \right\|
=
0,
\label{B*2B*0}
\end{eqnarray}
where 
${\bf B}^{*}(\alpha;0,x) = \left( Q(\alpha|x),0,\ldots,0\right)^T$.
\item
for all $(x',\theta_i)$ in $\mathcal{X} \times \Theta^1$, $i=1,2$,
\[
\left|
Q^{*}(x';\theta_1)
-
Q^{*}(x';\theta_2)
\right|
\leq
C
\underline{h}^{-p}
(
1+\underline{h}^{-1}
)
\left\|
\theta_1
-
\theta_2
\right\|
 .
\]

\item
There exists $C$ such that, for all $\theta$ in $\Theta^1$, all $x'$ in $\mathcal{X}$ and all $x$ in $\mathcal{X}_0$,
\[
f\left(
   Q^{*}
     (
      x';\theta
     )
     |
     x'
 \right)
K\left(\frac{x-x'}{h}\right)
\geq
C
K \left(\frac{x-x'}{h}\right)
.
\]
\end{enumerate}
\label{Biais1}
\end{lem}

\subsection{Proof of Theorem \protect{\ref{Bias}}}
Since $Q(\cdot|\cdot)$ is in $\mathcal{C}(L,s)$, the Taylor-Lagrange Formula and Assumption \ref{K} yield that there exists $t=t(h,x,z)$ in $(0,1)$ such that for $h$ small enough and all $(x,z)$ in $\mathcal{X}_0 \times \mathcal{K}$,
\begin{eqnarray}
Q(\alpha|x+hz)
& = &
\Sum_{0 \leq |{\bf v}| \leq \lfloor s \rfloor }
  \frac{
    b_{{\bf v}}(\alpha|x)
       }
       {{\bf v}!}
       (hz)^{{\bf v}}
+
\Sum_{|{\bf v}| = \lfloor s \rfloor} \frac{(hz)^{{\bf v}}}{{\bf v}!}
  \left(
    b_{{\bf v}}(\alpha|x + t h z)
      -
    b_{{\bf v}}(\alpha|x)
  \right)
\nonumber
\\
& = &
{\bf U}(z)^T {\bf H} {\bf b}(\alpha|x)
+
\epsilon (\theta,z).
\label{Rqthetaz}
\end{eqnarray}
In the equation above, $b_{{\bf v}} (\alpha|x)$ is the ${\bf v}$th partial derivatives of $Q(\alpha|x)$
with respect to $x$ and
${\bf b}(\alpha|x)  = \left(b_{{\bf v}} (\alpha|x), |{\bf v}| \leq \lfloor s \rfloor, 0,\ldots,0\right)^T \in \Rit^P$.
Since $Q(\cdot|\cdot) \in \mathcal{C} (L,s)$,
\begin{equation}
\lim_{\overline{h} \rightarrow 0}
\sup_{(\theta,z) \in \Theta^1 \times \mathcal{K}}
\left|
  \frac{
    \epsilon (\theta,z)
       }
       {h^s}
\right|
\leq
C L
.
\label{Rqbound}
\end{equation}

Let
\begin{eqnarray*}
I(\theta,z)
  =
\int_{0}^{1}
  f\left(
    Q(\alpha|x +hz)
    +
    t
     \left(
      {\bf U}(z)^T
      {\bf B}^{*}(\theta)
        -
      Q(\alpha|x+hz)
     \right)
     |x+hz
   \right)
dt
.
\end{eqnarray*}
Assumptions \ref{F}, \ref{K}, \ref{X}, $Q(\cdot|\cdot) \in \mathcal{C} (L,s)$ and (\ref{B*2B*0}) give
\begin{equation}
\lim_{\overline{h} \rightarrow 0}
\sup_{\left(\theta,z\right) \in \Theta^0\times \mathcal{K}}
\left|
  I(\theta,z)
  -
  f(Q(\alpha|x)|x)
\right|
=
0
.
\label{LimIthetaz}
\end{equation}
A Taylor expansion with integral remainder gives
\[
F\left(
  {\bf U}(z)^T
  {\bf B}^{*}(\theta)
  |x+hz
 \right)
-
  F\left(
  Q(\alpha|x+hz)
  |
  x+hz
   \right)
=
\left(
  {\bf U}(z)^T
  {\bf B}^{*}(\theta)
        -
  Q(\alpha|x+hz)
\right)
I(\theta,z)
.
\]
Substituting in the first-order condition (\ref{FOC}) yields
\begin{equation}
\int
  {\bf U}(z)
  \left(
    {\bf U}(z)^T 
    {\bf B}^{*}(\theta)
      -
    Q(\alpha|x+hz)
  \right)
  I(\theta,z)
  f(x+hz)
  K(z)
dz
=
0
 .
\label{Devtaylor}
\end{equation}
We show that the matrix 
$
\int
{\bf U}\left(z\right)
{\bf U}\left(z\right)^T
I\left(\theta,z\right)
f\left(x+hz\right)
K\left(z\right)
dz
$ 
has an inverse. Indeed, Assumptions \ref{K} and \ref{X}, (\ref{LimIthetaz}) and $\overline{h}$ small enough give that uniformly in $\theta$ in $\Theta^0$ and ${\bf A}$ in $\Rit^P$,
\begin{eqnarray*}
{\bf A}^T
\int
{\bf U}\left(z\right)
{\bf U}\left(z\right)^T
I\left(\theta,z\right)
f\left(x+hz\right)
K\left(z\right)
dz
{\bf A}
& = & 
\int
\left\|
{\bf U}\left(z\right)^T
{\bf A}
\right\|^2
I\left(\theta,z\right)
f\left(x+hz\right)
K\left(z\right)
dz
\\
& = &
\left(1+o(1)\right)
f\left(Q(\alpha|x)|x\right)
\int
\left\|
{\bf U}\left(z\right)^T
{\bf A}
\right\|^2
K\left(z\right)
dz
\\
& \geq & 
C 
\left\|
{\bf A}
\right\|^2,
\end{eqnarray*}
using the fact that
$
{\bf A}
\mapsto
\int
\left\|
{\bf U}\left(z\right)^T
{\bf A}
\right\|^2
K\left(z\right)
dz
$
is a square norm and norm equivalence over $\Rit^P$.
It follows that
$
\int
{\bf U}\left(z\right)
{\bf U}\left(z\right)^T
I\left(\theta,z\right)
f\left(x+hz\right)
K\left(z\right)
dz
$ 
is strictly positive definite and has an inverse which satisfies, for $n$ large enough
\begin{equation}
\sup_{\theta \in \Theta^0}
\left\|
\left[
\int
{\bf U}\left(z\right)
{\bf U}\left(z\right)^T
I\left(\theta,z\right)
f\left(x+hz\right)
K\left(z\right)
dz
\right]^{-1}
\right\|
< \infty.
\label{InvUUT}
\end{equation}
(\ref{Devtaylor}) and (\ref{Rqthetaz}) give
\[
{\bf H}{\bf b}^{*}(\theta)
=
{\bf H}{\bf b}(\alpha|x)
+
\left[\int {\bf U}(z){\bf U}(z)^T I(\theta,z) f(x+hz)K(z)dz\right]^{-1}
\int \epsilon (\theta,z) I(\theta,z) f(x+hz) {\bf U}(z) K(z)dz
 .
\]
It then follows from (\ref{Rqbound}) and (\ref{InvUUT}) that
\begin{eqnarray}
\lefteqn{
\left\|{\bf H}{\bf b}^{*}(\theta) -  {\bf H}{\bf b}(\alpha|x) \right\|
}
&&
\nonumber \\
& \leq &
\left\| \left[\int {\bf U}(z){\bf U}(z)^T I(\theta,z) f(x+hz)K(z)dz\right]^{-1} \right\|
\left\|\int \epsilon (\theta,z) I(\theta,z) f(x+hz) {\bf U}(z) K(z)dz \right\|
\nonumber \\
& \leq &
C L h^s
\label{Biascomplete}
\end{eqnarray}
uniformly in $\theta$ in 
$
\Theta^0$. This ends the proof of the Theorem and also establishes (\ref{B*highv}) since ${\bf b}(\alpha|x)  = \left(b_{{\bf v}} (\alpha|x), |{\bf v}| \leq \lfloor s \rfloor, 0,\ldots,0\right)^T$.
\eop

\subsection{Proof of Proposition \protect{\ref{Sharp}}}
Let $\varphi (t)= \exp(-t^2/2)/\sqrt{2 \pi}$, $\Phi (t) = \int_{-\infty}^{t} \varphi (u) du$ be the p.d.f  and c.d.f of the standard normal.
The regression model (\ref{Reg}) is such that
$$
F\left(y|x\right)=\Phi\left(y-m(x)\right),
\quad
f(x)=\ind\left(x \in \left[-1,1\right]\right).
$$
(\ref{B*2B*0}) gives that
$
\lim_{h\rightarrow 0}
\max_{z \in \mathcal{K}} 
\left|{\bf U} (z)^{T} {\bf B} (0.5;h,0) \right| = Q(0.5|0)=m(0)=0
$. 
Hence (\ref{Devtaylor}), (\ref{LimIthetaz}) and Assumption \ref{K} give
$$
\left(1+o(1) \right)
\varphi (0)
\int 
{\bf U} (z)
\left({\bf U} (z)^{T} {\bf B}(0.5;h,0)-m\left(hz \right)\right)
K(z) dz
=
0.
$$
Recall that 
$
{\bf U} (z)
=
\left(1,z\right)^T
$, so that the equation above gives
\begin{eqnarray*}
\left[
\begin{array}{l}
b_0 (0.5;h,0)
\\
h b_1 (0.5;h,0)
\end{array}
\right]
& = &
\left(
1+ o(1)
\right)
\left(
\int
{\bf U} (z) {\bf U}^{T} (z) K(z) dz
\right)^{-1}
\left[
\begin{array}{l}
h^{1/2} \int m(z) K(z) d z
\\
h^{1/2} \int |z|^{3/2} K(z) d z
\end{array}
\right]
\\
& = &
\left(
1+ o(1)
\right)
h^{1/2}
\left[
\begin{array}{l}
 \int m(z) K(z) d z
\\
\frac{\int |z|^{3/2} K(z) d z}{\int z^{2} K(z) d z}
\end{array}
\right]
.
\eop
\end{eqnarray*}

\subsection{Proof of Theorem \protect{\ref{Bahadur}}}

We first state some intermediary results. The two following propositions deals with the remainder term 
$
\mathbb{R}_n\left(\beta,\epsilon;\theta \right)
=
\sum_{i=1}^{n}
{\bf R}_i \left(\beta,\epsilon;\theta \right)
$ 
from (\ref{LL0}), where
\begin{eqnarray*}
\lefteqn{
{\bf R}_i \left(\beta,\epsilon;\theta \right)
}
&&
\\
& = &
\left\{
\ell_{\alpha}
\left(
  Y_i - Q^* (X_i ; \theta)
  -
  \frac{{\bf U} \left(\frac{X_i-x}{h}\right)^T\left(\beta+\epsilon\right)}{\left(n h^d\right)^{1/2}}
\right)
-
\ell_{\alpha}
\left(
  Y_i - Q^* (X_i ; \theta)
  -
  \frac{{\bf U} \left(\frac{X_i-x}{h}\right)^T\beta}{\left(n h^d\right)^{1/2}}
\right)\right\}
K
\left(\frac{X_i-x}{h}\right)
\\
&&
-
\frac{1}{\left(n h^d\right)^{1/2}}
{\bf S}_i (\theta)^{T} \epsilon
-
\frac{1}{2}
\epsilon^{T}
\left(
\frac{1}{n h^d}
{\bf J}_i (\theta)
\right) 
\left(\epsilon + 2 \beta \right).
\end{eqnarray*} 
Define also
\begin{eqnarray}
\label{Rho}
R_i \left(\beta,\epsilon;\theta \right)
&
=
& 
{\bf R}_i \left(\beta,\epsilon;\theta \right)
+
\frac{1}{2}
\epsilon^{T}
\left(
\frac{1}{n h^d}
{\bf J}_i (\theta)
\right) 
\left(\epsilon + 2 \beta \right)
\\
& = &
\Bigg\{
\ell_{\alpha}
\Bigg(
Y_i - Q^* (X_i ; \theta)
-
\frac{{\bf U} \left(\frac{X_i-x}{h}\right)^T\left(\beta+\epsilon\right)}{\left(n h^d\right)^{1/2}}
\Bigg)
\nonumber
\\
&&
-
\ell_{\alpha}
\Bigg(
  Y_i - Q^* (X_i ; \theta)
  -
  \frac{{\bf U} \left(\frac{X_i-x}{h}\right)^T\beta}{\left(n h^d\right)^{1/2}}
\Bigg)
\nonumber \\
&&
-
2
\left\{
\ind\left(Y_i \leq Q^* \left(X_i;\theta\right)\right)
-
\alpha
\right\}
\frac{{\bf U}
\left(\frac{X_i-x}{h}\right)^{T}
\epsilon}{\left(n h^d\right)^{1/2}}
\Bigg\}
K
\Big(\frac{X_i-x}{h}\Big),
\nonumber
\end{eqnarray}
\begin{eqnarray}
{\bf R}_i^{1} \left(\beta,\epsilon;\theta \right)
& = &
R_i \left(\beta,\epsilon;\theta \right)
-
\esp
\left[
R_i \left(\beta,\epsilon;\theta \right)
\left|X_i \right.
\right],
\label{R1}
\\
{\bf R}_i^{2} \left(\beta,\epsilon;\theta \right)
& = &
\esp
\left[
R_i \left(\beta,\epsilon;\theta \right)
\left|X_i \right.
\right]
-
\frac{1}{2}
\epsilon^{T}
\left(
\frac{1}{n h^d}
{\bf J}_i (\theta)
\right) 
\left(\epsilon + 2 \beta \right)
,
\label{R2}
\end{eqnarray}
which are such that 
$$
\mathbb{R}_n\left(\beta,\epsilon;\theta \right)
=
\mathbb{R}_n^{1} \left(\beta,\epsilon;\theta \right)
+
\mathbb{R}_n^{2} \left(\beta,\epsilon;\theta \right)
,
\quad
\mathbb{R}_n^j \left(\beta,\epsilon;\theta \right)
=
\sum_{i=1}^{n}
{\bf R}_i^{j} \left(\beta,\epsilon;\theta \right)
,
\; 
j=1,2.
$$
\begin{prop}
Consider two real numbers $t_{\beta},t_{\epsilon} > 0$ which may depend upon on $n$ with
$t_{\beta}\geq 1$, $t_{\epsilon}\geq 1/n$ and
$
\left(t_{\beta}+t_{\epsilon}\right)^{1/2}/t_{\epsilon}
\leq
O\left(\left(n \underline{h}^d\right)^{1/4}/ \log^{1/2} n\right)
$. Then, under Assumptions \ref{F}, \ref{K} and \ref{X} and for $n$ large enough, 
$$
\esp
\left[
\sup_{\left(\beta,\epsilon,\theta\right) \in \mathcal{B} (0,t_{\beta}) \times \mathcal{B} (0,t_{\epsilon}) \times \Theta^1 }
\left|
\mathbb{R}_n^{1} \left(\beta,\epsilon;\theta \right)
\right|
\right]
\leq
C
\frac{\log^{1/2}n}{\left(n \underline{h}^d\right)^{1/4}}
t_{\epsilon}\left(t_{\beta}+t_{\epsilon}\right)^{1/2}
.
$$ 
\label{R1order}
\end{prop}
\begin{prop}
Consider two real numbers $t_{\beta},t_{\epsilon} > 0$ which may depend upon on $n$ with
$t_{\beta}\geq 1$ and
$
t_{\beta}/t_{\epsilon}
=
O\left(n \underline{h}^d/ \log^{1/2} n\right)
$. Then, under Assumptions \ref{F}, \ref{K} and \ref{X} and for $n$ large enough, 
$$
\esp
\left[
\sup_{\left(\beta,\epsilon,\theta\right) \in \mathcal{B} (0,t_{\beta}) \times \mathcal{B} (0,t_{\epsilon}) \times \Theta^1 }
\left|
\mathbb{R}_n^{2} \left(\beta,\epsilon;\theta \right)
\right|
\right]
\leq
C
\frac{t_{\epsilon}\left(t_{\beta}+t_{\epsilon}\right)^{2}
}{
\left(n \underline{h}^d\right)^{1/2}}
.
$$ 
\label{R2order}
\end{prop}
The next lemma is used to bound the eigenvalues of $\sum^{n}_{i=1} \mathbf{J}_i (\theta)/(nh^d)$ from below. It implies in particular that all the $\beta_n (\theta)$ in (\ref{Betan}), $\theta$ in $\Theta^1$, are well defined  with a probability tending to 1. Let $\underline{\gamma}_n (\theta)$ be the smallest eigenvalue of the nonnegative symmetric matrix $\sum_{i=1}^{n} {\bf J}_i(\theta) / (n h)^d$.
\begin{lem}
Under Assumptions \ref{F}, \ref{K} and \ref{X},
$
\inf_{
  \theta \in \Theta^1
     } \underline{\gamma}_n (\theta)  \geq \underline{\gamma} + o_{\Prob}(1)
$ for some $\underline{\gamma}>0$.
\label{ControlJ}
\end{lem}
\noindent
Lemma \ref{ControlJ} together Lemma \ref{Ordrebetan} below gives 
$\sup_{\theta \in \Theta^1} \left\|\beta_n (\theta)\right\|=O_{\prob}
\left(\log ^{1/2} n \right)$.
\begin{lem}
Suppose that Assumptions \ref{F}, \ref{K} and \ref{X} are satisfied. Then 
$$
\sup_{\theta \in \Theta^1}
\left\|
\frac{1}{
\left(nh^d\right)^{1/2}
}
  \Sum_{i=1}^{n}
  {\bf S}_i (\theta)
\right\|
=
O_{\prob}
\left(\log ^{1/2} n \right).
$$
\label{Ordrebetan}
\end{lem}
The rest of the proof of Theorem \ref{Bahadur} is divided in two steps. In what follows
$$
t_n = t \frac{\log^{3/4} n}{ (n \underline{h}^d)^{1/4}},
\quad
t>0.
$$
Under Assumption \ref{K},
$(\log n) / (n \underline{h}^d ) =o(1)$
so that $t_n = o\left(\log^{1/2} n\right)$. In the sequel, $t_n$ will play the role of $t_{\epsilon}$ whereas $t_{\beta}$ will be chosen such that $t_{\beta} \asymp \log^{1/2} n$. Hence
\begin{eqnarray*}
\frac{\left( t_{\beta}+ t_{\epsilon} \right)^{1/2}}{t_{\epsilon}}
&\asymp&
\frac{(n \underline{h}^d)^{1/4}\log^{1/4} n}{t \log^{3/4} n}
=
\frac{1}{t} 
O\left(\frac{(n \underline{h}^d)^{1/4}}{\log^{1/2} n}\right),
\\
\frac{t_{\beta}}{t_{\epsilon}}
& \asymp &
\frac{(n \underline{h}^d)^{1/4} \log^{1/2} n}{t \log^{3/4} n}
=
O\left( \frac{n \underline{h}^d}{\log n} \right)^{1/4}
=
o
\left(
\frac{n \underline{h}^d}{\log n} \times \log^{1/2} n
\right)
=
o
\left(
\frac{n \underline{h}^d}{\log^{1/2} n}
\right).
\end{eqnarray*}
Hence these choices of $t_{\beta}$ and $t_{\epsilon}$ satisfy the conditions of Propositions \ref{R1order} and \ref{R2order} provided $t$ is chosen large enough.
\bigskip

{\it Step 1: order of 
$
\sup_{(\epsilon,\theta) \in \mathcal{B}(0,t_n) \times \Theta^1}
\left|
\mathbb{R}_n ( \beta_n (\theta),\epsilon;\theta)
\right|
$.} Consider $\eta>0$ arbitrarily small.
Let $\underline{\gamma}$ be as in Lemma \ref{ControlJ}. 
Since Lemmas \ref{ControlJ} and \ref{Ordrebetan} give
$\sup_{\theta \in \Theta^1} \left\|\beta_n (\theta)\right\|=O_{\prob}
\left(\log ^{1/2} n \right)$, there is a $C_{\eta}$ such that, for $n$ large enough,
\begin{eqnarray*}
\lefteqn{
\Prob
\left(
\sup_{(\epsilon,\theta) \in \mathcal{B}(0,t_n) \times \Theta^1}
\left|
\mathbb{R}_n ( \beta_n (\theta),\epsilon;\theta)
\right|
\geq
\frac{\underline{\gamma} t_n^2}{4}
\right)
}
& &
\\
& \leq &
\Prob
\left(
\sup_{(\epsilon,\theta) \in \mathcal{B}(0,t_n) \times \Theta^1}
\left|
\mathbb{R}_n ( \beta_n (\theta),\epsilon;\theta)
\right|
\geq
\frac{\underline{\gamma} t_n^2}{4}
,
\sup_{\theta \in \Theta^1} \left\|\beta_n (\theta)\right\|
\leq C_{\eta} \log^{1/2} n
\right)
\\
&&
+
\Prob
\left(
\sup_{\theta \in \Theta^1} \left\|\beta_n (\theta)\right\|
>
 C_{\eta} \log^{1/2} n
\right)
\\
& \leq &
\Prob
\left(
\sup_{(\beta,\epsilon,\theta) \in \mathcal{B}(0,C_{\eta}\log^{1/2}n )\times\mathcal{B}(0,t_n) \times \Theta^1}
\left|
\mathbb{R}_n ( \beta,\epsilon;\theta)
\right|
\geq
\frac{\underline{\gamma} t_n^2}{4}
\right)
+
\eta
.
\end{eqnarray*}
Propositions \ref{R1order} and \ref{R2order}, $\mathbb{R}_n = \mathbb{R}_n^1+\mathbb{R}_n^2$ and the Markov inequality give
\begin{eqnarray*}
\lefteqn{
\Prob
\left(
\sup_{(\beta,\epsilon,\theta) \in \mathcal{B}(0,C_{\eta}\log^{1/2}n )\times\mathcal{B}(0,t_n) \times \Theta^1}
\left|
\mathbb{R}_n ( \beta,\epsilon;\theta)
\right|
\geq
\frac{\underline{\gamma} t_n^2}{4}
\right)
}
& &
\\
& 
\leq
&  
\frac{C}{t_n^2}
\left(
\frac{
  t_n 
  \left( 
    C_{\eta}\log^{1/2} n + t_n     
  \right)^{1/2}
  \log^{1/2}n
  }
    {
      \left(
        n\underline{h}^d\right)^{1/4}
    }
+
\frac{
  t_n
  \left( 
    C_{\eta}\log^{1/2}n + t_n
  \right)^2
  }
  {
   \left( 
     n\underline{h}^d \right)^{1/2}}
   \right)
\\
& 
=
&
\frac{C}{t_n}
\frac{
  \log^{3/4} n
  }
  {
    (n \underline{h}^d)^{1/4}
  }
\left(
  \left( 
    C_{\eta} + \frac{t_n}{\log^{1/2} n}   
  \right)^{1/2}
+
  \left(
    \frac{\log n}
    {
      n \underline{h}^d
    } 
  \right)^{1/4}
  \left( 
    C_{\eta} + \frac{t_n}{\log^{1/2} n} 
  \right)^2
\right)
.
\end{eqnarray*}
The definition of $t_n$, $t_n=o\left(\log^{1/2} n \right)$ and Assumption \ref{K} give
\begin{equation}
\limsup_{n\rightarrow\infty}
\Prob
\left(
\sup_{(\epsilon,\theta) \in \mathcal{B}(0,t_n) \times \Theta^1}
\left|
\mathbb{R}_n ( \beta_n(\theta),\epsilon;\theta)
\right|
\geq
\frac{\underline{\gamma} t_n^2}{4}
\right)
=
\eta
+
O
\left(
\frac{C_{\eta}^{1/2}}{t}
\right)
\mbox{ \rm when $t \rightarrow\infty$.}
\label{Rorder}
\end{equation}

\bigskip

{\it Step 2: $\sup_{\theta \in \Theta^1} \left\|{\bf E}_n\left(\theta\right)\right\|$.}
Consider $\tau_n \geq t_n$ and $\epsilon = \tau_n {\bf e}$, $\|{\bf e}\|=1$ so that $\| \epsilon \|\geq t_n$. Since 
$\ell_{\alpha} (\cdot)$ is convex, $\epsilon \mapsto \mathbb{L}_n (\beta(\theta),\epsilon;\theta)$ is convex. This gives since $\mathbb{L}_n (\beta(\theta),0;\theta)=0$ and $\mathbb{L}_{n}=\mathbb{L}_{n}^0+\mathbb{R}_{n}$
\begin{eqnarray*}
  \frac{t_n}{\tau_n}
 \mathbb{L}_{n}
    \left(
      \beta_n(\theta) , \epsilon ; \theta
    \right)
 & = &
 \frac{t_n}{\tau_n}
 \mathbb{L}_{n}
    \left(
      \beta_n(\theta) , \epsilon ; \theta
    \right)
  +
  \left(
  1-\frac{t_n}{\tau_n}
  \right) 
  \mathbb{L}_{n}
    \left(
      \beta_n(\theta) , 0 ; \theta
    \right)
\\
& \geq &
\mathbb{L}_{n}
\left(
  \beta_n(\theta), \frac{t_n}{\tau_n} \epsilon
  ;
  \theta
\right)
=  
\mathbb{L}_{n}
\left(
  \beta_n(\theta), t_n {\bf e}
  ;
  \theta
\right)
\\
& \geq &
\mathbb{L}_{n}^0
\left(
  \beta_n(\theta), t_n {\bf e}
  ;
  \theta
\right)+\mathbb{R}_{n}\left(
  \beta_n(\theta), t_n {\bf e}
  ;
  \theta
\right)
.
\end{eqnarray*}
Hence ${\bf E}_n (\theta) = \arg \min_{\epsilon} \mathbb{L}_n (\beta_n (\theta),\epsilon;\theta)$ and the latter inequality give
\begin{eqnarray*}
\left\{\left\| {\bf E}_n(\theta)\right\| \geq t_n \right\}
& \subset &
\left\{
  \inf_{\epsilon ; \left\| \epsilon \right\| \geq t_n }
      \mathbb{L}_{n}(\beta_n(\theta),\epsilon;\theta)
    \leq 
    \inf_{\epsilon ; \left\| \epsilon \right\| < t_n }
      \mathbb{L}_{n}(\beta_n(\theta),\epsilon;\theta)
\right\}
\\
& \subset &
\left\{
  \inf_{\epsilon ; \left\| \epsilon \right\| \geq t_n }
      \mathbb{L}_{n}(\beta_n(\theta),\epsilon;\theta)
    \leq 
    \mathbb{L}_{n}(\beta_n(\theta),0;\theta)=0
\right\}
\\
& \subset &
\left\{
\inf_{\mathbf{e} ; \left\| \mathbf{e} \right\| = 1 }
\left[
\mathbb{L}_{n}^0
\left(
\beta_n(\theta), t_n {\bf e}
;
\theta
\right)
+
\mathbb{R}_{n}
\left(
\beta_n(\theta), t_n {\bf e}
 ;
\theta
\right)
\right]
\leq 
0
\right\}
\\
& \subset &
\left\{
\inf_{ \left\| \epsilon \right\| =t_n}
\mathbb{L}_{n}^{0} \left(\beta_n(\theta),\epsilon;\theta\right)
-
\sup_{\left\|\epsilon\right\| = t_n}
\left|
\mathbb{R}_{n}
\left(
\beta_n(\theta), \epsilon
;
\theta
\right)
\right|
\leq
0
\right\}
.
\end{eqnarray*}
Since
\begin{eqnarray*}
\left\{
  \Sup_{\theta \in \Theta^1}
    \left\|
      {\bf E}_n(\theta)
    \right\|
  \geq
    t_n
\right\}
=
\bigcup_{\theta \in \Theta^1}
  \left\{
    \left\|
      {\bf E}_n(\theta)
    \right\|
  \geq
  t_n
\right\}
,
\end{eqnarray*}
this gives
\begin{eqnarray}
\left\{
  \Sup_{\theta \in \Theta^1}
    \left\|
      {\bf E}_n(\theta)
    \right\|
  \geq
  t_n
\right\}
& \subset &
\bigcup_{\theta \in \Theta^1}
\left\{
\inf_{ \left\| \epsilon \right\| =t_n}
\mathbb{L}_{n}^{0} \left(\beta_n(\theta),\epsilon;\theta\right)
-
\sup_{\left\|\epsilon\right\| = t_n}
  \left|
    \mathbb{R}_{n}\left(
  \beta_n(\theta), \epsilon
  ;
  \theta
\right)
   \right|
\leq
0
\right\}
\nonumber
\\
& \subset &
\left\{
\inf_{\theta \in \Theta^1}
\inf_{ \left\| \epsilon \right\| =t_n}
\mathbb{L}_{n}^{0} \left(\beta_n(\theta),\epsilon;\theta\right)
\leq
\sup_{(\epsilon,\theta) \in \mathcal{B}(0,t_n) \times \Theta^1}
  \left|
    \mathbb{R}_{n}\left(
  \beta_n(\theta), \epsilon
  ;
  \theta
\right)
\right|
\right\}
.
\label{Inclusion}
\end{eqnarray}
Consider first
$
\inf_{\theta \in \Theta^1}
\inf_{ \left\| \epsilon \right\| =t_n}
    \mathbb{L}_{n}^{0} \left(\beta_n(\theta),\epsilon;\theta\right)
$. The definition (\ref{LL0}) of $\mathbb{L}_{n}^{0}$ gives, for any $\epsilon$ with 
$\left\| \epsilon \right\| =t_n$,
$$
\mathbb{L}_{n}^{0} \left(\beta_n(\theta),\epsilon;\theta\right)
=
\frac{1}{2}
\epsilon^T
     \left(
      \frac{1}{nh^d}
      \Sum_{i=1}^{n} {\bf J}_i(\theta)
    \right)
\epsilon
\geq
\frac{1}{2}
\underline{\gamma}_n (\theta)
t_n^2.
$$
Hence (\ref{Inclusion}), Lemma \ref{ControlJ} and (\ref{Rorder}) give
\begin{eqnarray*}
\limsup_{n\rightarrow\infty}
\Prob
\left(
\Sup_{\theta \in \Theta^1}
    \left\|
      {\bf E}_n(\theta)
    \right\|
  \geq
  t_n
\right)
& \leq &
\limsup_{n\rightarrow\infty}
\Prob
\left(
\sup_{(\epsilon,\theta) \in \mathcal{B}(0,t_n) \times \Theta^1}  \left|
    \mathbb{R}_{n}\left(
  \beta_n(\theta), \epsilon
  ;
  \theta
\right)
\right|
\geq
\frac{\underline{\gamma}_n (\theta) t_n^2}{2}
\right)
\\
& \leq &
\limsup_{n\rightarrow\infty}
\Prob
\left(
\sup_{(\epsilon,\theta) \in \mathcal{B}(0,t_n) \times \Theta^1}  
\left|
    \mathbb{R}_{n}\left(
  \beta_n(\theta), \epsilon
  ;
  \theta
\right)
\right|
\geq
\frac{\underline{\gamma}t_n^2}{4}
\right)
\\
& = &
\eta
+
O
\left(
\frac{C_{\eta}^{1/2}}{t}
\right)
\mbox{ \rm when $t\rightarrow\infty$.}
\end{eqnarray*}
Since the latter can be made arbitrarily small by taking $\eta$ arbitrarily small and then $t$ large enough,
the Theorem is proved. \eop

\subsection{Proof of Corollary \protect{\ref{Optimal}}} Part (i) follows from Theorems \ref{Bias} and \ref{Bahadur} and the triangular inequality, together with
$$
\left(
\int_{\mathcal{X}_0}
\left\|
\beta_n (\alpha;h,x)
\right\|^{m}
dx
\right)^{1/m}
=
O_{\prob} (1).
$$
We now prove the latter. Lemma \ref{ControlJ}
and the Hölder inequality give, since $\mathcal{X}_0$ is compact,
$$
\left(
\int_{\mathcal{X}_0}
\left\|
\beta_n (\alpha;h,x)
\right\|^{m}
dx
\right)^{1/m}
=
O_{\prob}
\left(
\int_{\mathcal{X}_0}
\left\|
\frac{1}{(nh^d)^{1/2}}
\sum_{i=1}^n
\mathbf{S}_i (\alpha;h,x)
\right\|^{2[m]+2}
dx
\right)^{1/(2[m]+2)}
.
$$
Since $\esp [ {\bf S}_i (\theta) ] = 0$, the  Marcinkiewicz-Zygmund inequality (see Chow and Teicher, 2003), (\ref{Score}) and $h^d \geq C (\log n)/n$ give
\begin{eqnarray*}
\lefteqn{
\esp^{1/(2[m]+2)}
\left[
\left\|
\frac{1}{(nh^d)^{1/2}}
\sum_{i=1}^n
\mathbf{S}_i (\theta)
\right\|^{2[m]+2}
\right]
\leq 
C
\esp^{1/(2[m]+2)}
\left[
\left(
\frac{1}{nh^d}
\sum_{i=1}^n
\left\|
\mathbf{S}_i (\theta)
\right\|^2
\right)^{[m]+1}
\right]
}
\\
& \leq &
C
\left(
\frac{1}{(nh^d)^{[m]+1}}
\sum_{i_1,\ldots,i_{[m]+1}=1}^n
\esp
\left[
\ind
\left(
\frac{X_{i_1}-x}{h}
\in 
\mathcal{K}
\right)
\times
\cdots
\times
\ind
\left(
\frac{X_{i_{[m]+1}}-x}{h}
\in 
\mathcal{K}
\right)
\right]
\right)^{1/2}
=
O(1),
\end{eqnarray*}
uniformly in $x$.
Part (ii) similarly follows from Lemmas \ref{ControlJ} and \ref{Ordrebetan} which gives
$
\sup_{\theta \in \Theta^1}
\left\|
\beta_n (\theta)
\right\|
=
O_{\prob} \left( \log^{1/2} n \right)
$.
\eop

\subsection{Proof of Proposition \protect{\ref{Randh}}}
Let $\underline{h} = h_n/C$ and $\overline{h} = C h_n$. The condition on $h_n$ ensures that $\underline{h}$ and $\overline{h}$ satisfy Assumption \ref{K} for all $C>1$.
Recall that
Lemma \ref{ControlJ} together Lemma \ref{Ordrebetan} gives 
$\sup_{\theta \in \Theta^1} \left\|\beta_n (\theta)\right\|=O_{\prob}
\left(\log ^{1/2} n \right)$.
Hence (\ref{En}), Theorems \ref{Bias} and \ref{Bahadur} give, for all $C>1$,
\begin{eqnarray*}
\sup_{
(\alpha,x,h) \in [\underline{\alpha},\overline{\alpha}] \times \mathcal{X}_0 \times [\underline{h},\overline{h}]
}
\left|
\widehat{b}_{\bf v} (\alpha; \widehat{h},x)
-
b_{\bf v} (\alpha | x)
\right|
& = &
\underline{h}^{-|{\bf v}|}
O_{\prob}
\left(
\underline{h}^s
+
\left(
\frac{\log n}{n \underline{h}^d}
\right)^{1/2}
\right)
\\
& = &
h_n^{-|{\bf v}|}
O_{\prob}
\left(
h_n^s
+
\left(
\frac{\log n}{n h_n^d}
\right)^{1/2}
\right)
\; .  
\end{eqnarray*}
This ends the proof of the Proposition since 
$\liminf_{n\rightarrow\infty} \Prob\left(\widehat{h}_n \in [\underline{h},\overline{h}]\right)
$
can be made arbitrarily close to 1 by increasing $C$. \eop

\subsection{Proof of Proposition \protect{\ref{Quantdens}}}
Substituting (\ref{En}) in (\ref{Hatq}) yields
\begin{eqnarray*}
\widehat{q} (\alpha|x) - q (\alpha|x)
& = &
\frac{1}{h_q}
\int Q(\alpha+h_qt|x)
d K_q (t)
- 
q (\alpha|x)
\\
&&
+
\frac{1}{h_q}
\int 
\left(
Q^*(\alpha+h_qt|x)
-
Q(\alpha+h_qt|x)
\right)
d K_q (t)
\\
& &
+
\int
\frac{
\mathbf{e}_0^{T} 
\beta_n
\left(
\alpha+h_qt
; h,x
\right)}{h_q \left( n h^d \right)^{1/2}}
d K_q (t)
+
\int
\frac{\mathbf{e}_0^{T} 
\mathbf{E}_n
\left(
\alpha+h_qt
; h,x
\right)}{h_q\left( n  h^d \right)^{1/2}}
d K_q (t)
.
\end{eqnarray*}
Theorems \ref{Bias} and \ref{Bahadur} with $h=O(h_q)$ and $h_q\rightarrow0$ give
\begin{eqnarray*}
\frac{1}{h_q}
\int 
\left(
Q^*(\alpha+h_qt|x)
-
Q(\alpha+h_qt|x)
\right)
d K_q (t)
& = &
O
\left(
\frac{h^{s+1}}{h_q}
\int
\left|d K_q (t)\right|
\right)
=
O (h_q^s),
\\
\int 
\frac{
\mathbf{E}_n
\left(
\alpha+h_qt
; h,x
\right)}{h_q \left( n  h^d \right)^{1/2}}
d K_q (t)
& = &
\frac{\log^{3/4} n}{\left( nh^d h_q^2\right)^{1/4}}
O_{\prob}
\left(
\frac{1}{\left( n h_q h^d \right)^{1/2}}
\right)
.
\end{eqnarray*}
Hence it remains to show that
\begin{eqnarray}
\frac{1}{h_q}
\int Q(\alpha+h_qt|x)
d K_q (t)
- 
q (\alpha|x)
& = &
O\left(h_q^s\right),
\label{Quantdens1}
\\
\frac{1}{h_q^{1/2}}
\int 
\beta_n
\left(
\alpha+h_qt
; h,x
\right)
d K_q (t)
& = &
O_{\prob} (1)
.
\label{Quantdens2}
\end{eqnarray}
The two next steps establish these two equalities.

\medskip

{\it Step 1: proof of (\ref{Quantdens1})}. Let $q^{(j)}(\alpha|x)= \partial^j q(\alpha|x) / \partial x^j$.
Since $Q(\alpha|x) \in \mathcal{C} (L,s+1)$, the Taylor-Lagrange Formula gives, for some $\omega$ in $[0,1]$,
$$
Q(\alpha + h_q t|x)
-
Q(\alpha|x)
=
\sum_{j=0}^{\left\lfloor s\right\rfloor}
\frac{q^{(j)}(\alpha|x)}{(j+1)!}
\left(h_q t\right)^j
+
\frac{
q^{(\left\lfloor s\right\rfloor)}(\alpha+\omega h_q t|x)
-
q^{(\left\lfloor s\right\rfloor)}(\alpha|x)}{(\left\lfloor s\right\rfloor+1)!}
\left(h_q t\right)^{\left\lfloor s\right\rfloor}
.
$$
The definition of the smoothness class $\mathcal{C}_q (L,s)$ gives
$$
\left|
q^{(\left\lfloor s\right\rfloor)}(\alpha+\omega h_q t|x)
-
q^{(\left\lfloor s\right\rfloor)}(\alpha|x)
\right|
\leq
L \left|h_q t\right|^{s-\left\lfloor s\right\rfloor}
.
$$
Hence, since the support of $K_q (\cdot)$ is compact, $\int\left|dK_q (t)\right|<\infty$ and
$\int dK_q (t) = 0$, $\int t dK_q (t) = 1$, $\int t^2 dK_q (t) = \cdots = \int t^{\left\lfloor s\right\rfloor} dK_q (t)=0$,
\begin{eqnarray*}
\frac{1}{h_q}
\int Q(\alpha+h_qt|x)
d K_q (t)
& = &
\frac{Q(\alpha|x)}{h_q}
\int 
d K_q (t)
+
q(\alpha|x)
\int
t 
d K_q (t)
+
\frac{h_q q^{(1)}(\alpha|x)}{2}
\int
t^2 
d K_q (t)
\\
&&
+\cdots+
\frac{h_q^{\left\lfloor s\right\rfloor} q^{\left(\left\lfloor s\right\rfloor\right)}(\alpha|x)}{\left(\left\lfloor s\right\rfloor+1\right)}
\int
t^{\left\lfloor s\right\rfloor} 
d K_q (t)
+
O(h^s)
\\
& = &
q(\alpha|x) + O(h^s).
\end{eqnarray*}

\medskip

{\it Step 2: proof of (\ref{Quantdens2})}. Let $\theta_t = \left(\alpha+h_qt,h,x\right)$, $\theta=\theta_0$.
Since $\int dK_q (t) = 0$, (\ref{Betan}) gives
\begin{eqnarray}
\lefteqn{
\frac{1}{h_q^{1/2}}
\int 
\beta_n
\left(
\theta_t
\right)
d K_q (t)
= 
\frac{1}{h_q^{1/2}}
\int 
\left(
\beta_n
\left(
\theta_t
\right)
-
\beta_n \left(\theta\right)
\right)
d K_q (t)
}
\nonumber
&&
\\
& = &
\frac{1}{h_q^{1/2}}
\int 
\left\{
\left(
\frac{1}{n h^d}
\sum_{i=1}^n
\mathbf{J}_i
\left(
\theta
\right)
\right)^{-1}
-
\left(
\frac{1}{n h^d}
\sum_{i=1}^n
\mathbf{J}_i
\left(
\theta_t
\right)
\right)^{-1}
\right\}
\frac{1}{\left(n h^d\right)^{1/2}}
\sum_{i=1}^n
\mathbf{S}_i (\theta)
d K_q (t)
\label{DeltJ}
\\
&&
+
\frac{1}{h_q^{1/2}}
\int
\left(
\frac{1}{n h^d}
\sum_{i=1}^n
\mathbf{J}_i
\left(
\theta_t
\right)
\right)^{-1}
\frac{1}{\left(n h^d\right)^{1/2}}
\sum_{i=1}^n
\left\{
\mathbf{S}_i
\left(
\theta
\right)
-
\mathbf{S}_i
\left(
\theta_t
\right)
\right\}
d K_q (t).
\label{DeltS}
\end{eqnarray}
Since $\mathbf{A} \mapsto \mathbf{A}^{-1}$ is Lipshitz over the set of semi-definite positive matrices $\mathbf{A}$ with smallest eigenvalue bounded from below by $\underline{\gamma}$, Lemmas \ref{ControlJ} and A.3, (\ref{J}) and Assumption \ref{F} yield that (\ref{DeltJ}) satisfies
\begin{eqnarray*}
\lefteqn{
\left\|
\frac{1}{h_q^{1/2}}
\int 
\left\{
\left(
\frac{1}{n h^d}
\sum_{i=1}^n
\mathbf{J}_i
\left(
\theta
\right)
\right)^{-1}
-
\left(
\frac{1}{n h^d}
\sum_{i=1}^n
\mathbf{J}_i
\left(
\theta_t
\right)
\right)^{-1}
\right\}
\frac{1}{\left(n h^d\right)^{1/2}}
\sum_{i=1}^n
\mathbf{S}_i (\theta)
d K_q (t)
\right\|
}
&&
\\
& \leq &
\frac{O_{\prob}(1)}{h_q^{1/2}}
\int
\left\|
\frac{1}{n h^d}
\sum_{i=1}^n
\left\{
\mathbf{J}_i
\left(
\theta_t
\right)
-
\mathbf{J}_i
\left(
\theta
\right)
\right\}
\right\|
\left\|
\frac{1}{\left(n h^d\right)^{1/2}}
\sum_{i=1}^n
\mathbf{S}_i (\theta)
\right\|
\left|
d K_q (t)
\right|
\\
& \leq &
\frac{O_{\prob}(\log n)^{1/2}}{h_q^{1/2}}
\int
\frac{1}{n h^d}
\sum_{i=1}^n
\left|
Q^*\left(X_i;\theta_t\right)
-
Q^*\left(X_i;\theta\right)
\right|
\ind
\left(
\frac{X_i-x}{h}
\in 
\mathcal{K}
\right)
\left|
d
K_q (t)
\right|.
\end{eqnarray*}
The definition (\ref{Q*}) of $Q^*(X;\theta)$ and (\ref{Biascomplete}) give, since $Q(\alpha|x) \in \mathcal{C} (L,s+1)$ and because the support of $K_q $ is compact,
\begin{eqnarray*}
\lefteqn{
\frac{1}{h_q^{1/2}}
\int
\frac{1}{n h^d}
\sum_{i=1}^n
\left|
Q^*\left(X_i;\theta_t\right)
-
Q^*\left(X_i;\theta\right)
\right|
\ind
\left(
\frac{X_i-x}{h}
\in 
\mathcal{K}
\right)
\left|
d
K_q (t)
\right|
}
&&
\\
&=&
\frac{1}{h_q^{1/2}}
\int
\frac{1}{n h^d}
\sum_{i=1}^n
\left|
\mathbf{U} \left(\frac{X_i-x}{h}\right)^T
\left(
\mathbf{H}
\mathbf{b}^*(\theta_t)
-
\mathbf{H}
\mathbf{b}^*(\theta)
\right)
\right|
\ind
\left(
\frac{X_i-x}{h}
\in 
\mathcal{K}
\right)
\left|
d
K_q (t)
\right|
\\
& \leq &
C
\frac{1}{n h^d}
\sum_{i=1}^n
\ind
\left(
\frac{X_i-x}{h}
\in 
\mathcal{K}
\right)
\frac{1}{h_q^{1/2}}
\int
\left\|
\mathbf{H}
\mathbf{b}^*(\theta_t)
-
\mathbf{H}
\mathbf{b}^*(\theta)
\right\|
\left|
d
K_q (t)
\right|
\\
&\leq&
O_{\prob} (1)
\frac{1}{h_q^{1/2}}
\left(
\int
\left\|
\mathbf{H}
\mathbf{b} (\alpha+h_q t|x)
-
\mathbf{H}
\mathbf{b} (\alpha|x)
\right\|
\left|
d
K_q (t)
\right|
+
O(h^{s+1})
\right)
\\
& \leq &
O_{\prob} (1)
\frac{1}{h_q^{1/2}}
\left(
\int
\left|
Q\left(\alpha+h_qt|x\right)
-
Q\left(\alpha|x\right)
\right|
\left|
d
K_q (t)
\right|
+
O
\left(
h^{s+1}
+
h
\right)
\right)
=O_{\prob} \left(h_q^{1/2}\right).
\end{eqnarray*}
This gives that the item in (\ref{DeltJ}) is $O_{\prob} \left(h_q^{1/2}\right)=o_{\prob} (1)$.

For (\ref{DeltS}), Lemma \ref{ControlJ}, $\esp [\mathbf{S}_i (\theta_t)]=0$, (\ref{Score}), $Q^* (X;\theta_t)=Q^*(X;\theta)+O(h_q)$ uniformly with respect to $t$ in the support of $K_q$ and $X \in x + h \mathcal{K}$ (as easily seen arguing as in the equation above) and Assumptions
\ref{F}, \ref{X} give
\begin{eqnarray*}
\lefteqn{
\left\|
\frac{1}{h_q^{1/2}}
\int
\left(
\frac{1}{n h^d}
\sum_{i=1}^n
\mathbf{J}_i
\left(
\theta_t
\right)
\right)^{-1}
\frac{1}{\left(n h^d\right)^{1/2}}
\sum_{i=1}^n
\left\{
\mathbf{S}_i
\left(
\theta
\right)
-
\mathbf{S}_i
\left(
\theta_t
\right)
\right\}
d K_q (t)
\right\|
}
&&
\\
&\leq &
\frac{O_{\prob} (1)}{h_q^{1/2}}
\int
\left\|
\frac{1}{\left(n h^d\right)^{1/2}}
\sum_{i=1}^n
\left\{
\mathbf{S}_i
\left(
\theta
\right)
-
\mathbf{S}_i
\left(
\theta_t
\right)
\right\}
\right\|
\left|
d K_q (t)
\right|
\\
& = &
\frac{O_{\prob} (1)}{h_q^{1/2}}
\esp
\left[
\int
\left\|
\frac{1}{\left(n h^d\right)^{1/2}}
\sum_{i=1}^n
\left\{
\mathbf{S}_i
\left(
\theta
\right)
-
\mathbf{S}_i
\left(
\theta_t
\right)
\right\}
\right\|
\left|
d K_q (t)
\right|
\right]
\\
& \leq &
\frac{O_{\prob} (1)}{h_q^{1/2}}
\int
\esp^{1/2}
\left[
\left\|
\frac{1}{\left(n h^d\right)^{1/2}}
\sum_{i=1}^n
\left\{
\mathbf{S}_i
\left(
\theta
\right)
-
\mathbf{S}_i
\left(
\theta_t
\right)
\right\}
\right\|^2
\right]
\left|
d K_q (t)
\right|
\\
& = &
\frac{O_{\prob} (1)}{h_q^{1/2}}
\int
\Var^{1/2}
\left(
\frac{1}{h^{d/2}}
\left\{
\mathbf{S}
\left(
\theta
\right)
-
\mathbf{S}
\left(
\theta_t
\right)
\right\}
\right)
\left|
d K_q (t)
\right|
\\
& = &
O_{\prob} (1)
\int
\left[
\int
\left(
\frac{1}{h_q}
\int_{Q^*(x+hz;\theta)-Ch_q}^{Q^*(x+hz;\theta)+Ch_q}
f(y|x+hz)
dy
\right)
\ind \left( z \in \mathcal{K} \right)
f(x+hz)
dz
\right]^{1/2}
\left|
d K_q (t)
\right|
\\
& = &
O_{\prob} (1). \eop
\end{eqnarray*}

\setcounter{lem}{0} 
\setcounter{prop}{0}
\setcounter{thm}{0} 
\renewcommand{\thelem}{B.\arabic{lem}} 
\renewcommand{\theprop}{B.\arabic{prop}} 
\setcounter{equation}{0} \setcounter{subsection}{0} 
\renewcommand{\theequation}{B.\arabic{equation}} 
\renewcommand{\thesubsection}{B.\arabic{subsection}}
\section*{Appendix B: Proofs of intermediary results}

\subsection{Proof of Lemma \protect{\ref{Biais1}}}
Recall
\[
{\bf U}(X-x)^T {\bf b} = {\bf U}(X-x)^T {\bf H}^{-1} {\bf B} = {\bf U}\left(\frac{X-x}{h}\right)^T {\bf B}
\]
and define 
\[
\widetilde{\mathcal{L}}({\bf B};\theta)
=
\mathcal{L} \left( {\bf b};\theta \right)
=
  \frac{1}{h^d}
  \esp
  \left[
       \left\{
         \ell_{\alpha}
           (Y - {\bf U}((X-x)/h)^T {\bf B})
           -
         \ell_{\alpha}
           (Y)
        \right\}
       K_{h}(X-x)
  \right]
.
\]
The change of variable $x_1 = x+hz$ gives
\begin{eqnarray}
\widetilde{\mathcal{L}}({\bf B};\theta)
 & = &
  \frac{1}{h^d}
   \int
     \left[
       \int
         \left(
           \ell_{\alpha}(y - {\bf U}\left(\frac{x_1-x}{h}\right)^T {\bf B})
             -
           \ell_{\alpha}(y)
          \right)
          f(y|x_1)
          dy
        \right]
       f(x_1)
       K\left(\frac{x_1-x}{h}\right)
       dx_1
\nonumber
\\
& = &
   \int
     \left[
       \int
         \left(
           \ell_{\alpha}(y - {\bf U}\left(z\right)^T {\bf B})
             -
           \ell_{\alpha}(y)
          \right)
          f(y|x +hz)
          dy
        \right]
       f(x+hz)
       K\left(z\right)
       dz
,
\label{Gtilde}
\end{eqnarray}
showing that $\widetilde{\mathcal{L}} \left({\bf B};\theta\right)$ is also defined for $h=0$.

{\it Proof of (i).} It is sufficient to show that ${\bf B}^{*}(\theta) = \arg\min_{{\bf B} \in \Rit^P} \widetilde{\mathcal{L}} ({\bf B};\theta)$ exists and is unique. Note that $\mathbf{B} \mapsto \widetilde{\mathcal{L}} ({\bf B};\theta)$ is convex  by (\ref{Gtilde}) because $\ell_{\alpha} (\cdot)$ is convex. 
Since $\lim_{|t| \rightarrow +\infty} \ell_{\alpha}(t) = +\infty$
and
${\bf U}(z)^T {\bf B}$ diverges almost everywhere  when $\left\|{\bf B}\right\|$ diverges, 
(\ref{Gtilde}) gives that $\lim_{\left\|{\bf B}\right\| \rightarrow + \infty} \widetilde{\mathcal{L}}({\bf B};\theta) = +\infty$. Hence $\widetilde{\mathcal{L}}({\bf B};\theta)$ has a minimum. We  show that this minimum is unique by showing that $\mathbf{B} \mapsto \widetilde{\mathcal{L}}({\bf B};\theta)$ is strictly convex for all $\theta$ in $\Theta^0$. We compute the first and second ${\bf B}$-derivatives 
 of $\widetilde{\mathcal{L}}({\bf B};\theta)$. Equation (\ref{Lalpha}) gives that for almost all ${\bf B}$,
\[
\frac{
  \partial
    \ell_{\alpha}
      \left(
        y - {\bf U}(z)^T {\bf B}
      \right)
      }
      {
        \partial {\bf B}^T
      }
=
2
  \left(
    \ind
      \left(
        y \leq {\bf U}(z)^T {\bf B}
      \right)
     -
     \alpha
   \right)
{\bf U}(z)
\]
which is bounded for $z$ in the compact $\mathcal{K}$. Assumptions \ref{F}, \ref{K} and \ref{X}, the Lebesgue Dominated Convergence Theorem and (\ref{Gtilde})  yield that
\begin{eqnarray}
\widetilde{\mathcal{L}}^{(1)}({\bf B};\theta)
& = &
\frac{\partial \widetilde{\mathcal{L}}({\bf B};\theta) }{\partial {\bf B}^T}
=
2
 \int
   \left(
     \int
       \left(
         \ind
           \left(
             y \leq {\bf U}(z)^T {\bf B}
           \right)
             -
           \alpha
       \right)
       f(y|x+hz)
       dy
   \right)
   f(x+hz)
   {\bf U}(z)
   K(z)
dz
\nonumber
\\
& = &
  2
   \int
     F\left({\bf U}\left(z\right)^T {\bf B} | x+hz \right)
     f(x+hz)
     {\bf U}(z)
     K(z)
   dz
 -2
 \alpha
 \int
   f\left(x+hz\right)
   {\bf U}(z)
   K\left(z\right)
 dz
 .
\label{Gtilde1}
\end{eqnarray}
Applying again the Dominated Convergence Theorem yields that
\begin{equation}
\widetilde{\mathcal{L}}^{(2)} ({\bf B};\theta)
=
\frac{\partial^2 \widetilde{\mathcal{L}} ({\bf B};\theta)}{\partial {\bf B}^T \partial {\bf B}}
=
2
  \int
    f({\bf U}(z)^T{\bf B}|x+hz)f(x+hz){\bf U}(z){\bf U}(z)^T K(z)dz
 .
\label{Gtilde2}
\end{equation}
For all ${\bf A} \neq 0$ in $\Rit^P$, (\ref{Gtilde2}),  Assumptions \ref{F}, \ref{K} \ref{X} and $x \in \mathcal{X}_0$ give
\begin{eqnarray}
{\bf A}^T
\widetilde{\mathcal{L}}^{(2)}
  ({\bf B};\theta)
{\bf A}
& = &
2
  \int
    f(
       {\bf U}(z)^T
       {\bf B}
       |
       x+hz
      )
    f(
      x+hz
     )
     {\bf A}^T
     {\bf U}(z)
     {\bf U}(z)^T
     {\bf A}
     K(z)
dz
\nonumber
\\
& = &
2
\int
  f\left({\bf U}(z)^T {\bf B} |x+hz\right)
  f\left(x+hz\right)
  \left\|
    {\bf U}\left(z\right)^T
    {\bf A}
  \right\|^2
  K\left(z\right)
dz
>
0.
\label{Formquad}
\end{eqnarray}
Hence $\widetilde{\mathcal{L}}^{(2)}(\cdot;\theta)$ is a positive definite symmetric matrix for all $\theta$ in $\Theta^0$ and ${\bf B}$ in $\Rit^P$ so that the strictly convex function $\widetilde{\mathcal{L}}({\bf B};\theta)$ achieves it minimum for a unique ${\bf B}^{*}(\theta)$.

{\it Proof of (ii).} 
Consider a fixed $\overline{h}$ to be chosen small enough, and let $\widetilde{\Theta}^0$ be the corresponding $\Theta^0$, which is compact. The proof of (i) yields that ${\bf B}^{*}(\theta)$ is unique for all $\theta$ in $\widetilde{\Theta}^0$ and is the unique solution of the first-order condition
$\widetilde{\mathcal{L}}^{(1)}({\bf B};\theta)=0$, that is
\begin{equation}
\int
   F\left({\bf U}\left(z\right)^T {\bf B}| x+hz \right)
   f(x+hz)
   {\bf U}(z)
   K(z)
 dz
=
 \alpha
  \int
    f\left(x+hz\right)
    {\bf U}(z)
    K\left(z\right)
  dz
\label{CPOB}
,
\end{equation}
see (\ref{Gtilde1}), so that (\ref{FOC}) is proved.
In particular, ${\bf B}^{*}(\alpha;0,x)$ is the unique solution of $\widetilde{\mathcal{L}}^{(1)}({\bf B};\alpha,0,x) = 0$. If $h=0$, the first order condition (\ref{FOC}) is equivalent to
\[
\int
  F({\bf U}(z)^T {\bf B}^{*}(\alpha;0,x)|x)
  {\bf U}(z)
  K(z)
dz
=
\alpha
\int
  {\bf U}(z)
  K(z)
dz
.
\]
Let ${\bf B}_0^T(\alpha|x) = (Q(\alpha|x),0,\ldots,0)$ in $\Rit^P$. Since ${\bf U}(z)^T {\bf B}_0(\alpha|x) = Q(\alpha|x)$, ${\bf B}_0(\alpha|x)$ satisfies the first-order condition equation above.
Hence ${\bf B}^{*}(\alpha;0,x) = {\bf B}_0(\alpha|x)$ by uniqueness.

We now show that ${\bf B}^{*}(\theta)$ is continuously differentiable in $\theta$ over $\widetilde{\Theta}^0$ and give bounds for ${\bf B}^{*}(\theta)$, $\partial \widetilde{\mathcal{L}}^{(1)} ({\bf B}^{*}(\theta);\theta) / \partial \theta^T$ and $\widetilde{\mathcal{L}}^{(2)} ({\bf B}^{*}(\theta);\theta)$.
As shown above, ${\bf B} \mapsto \widetilde{\mathcal{L}}^{(1)}({\bf B};\theta)$ is continuously differentiable and $\widetilde{\mathcal{L}}^{(2)}({\bf B};\theta)$ is a symmetric positive definite matrix for all ${\bf B}$ in $\Rit^P$ and so has an inverse.
Assumptions \ref{F}, \ref{K} and \ref{X} yield that $F({\bf U}(z)^T {\bf B} | x+hz )$ and $f(x+hz)$ are bounded and have bounded $\theta$-partial derivatives over $\widetilde{\Theta}^0$ provided  $\overline{h}$ is small enough.
Hence the Dominated Convergence Theorem  and (\ref{Gtilde1}) yield that $\widetilde{\mathcal{L}}^{(1)}({\bf B};\theta)$ is continuously differentiable in $\theta$ over $\widetilde{\Theta}^0$.
Then the Implicit Function Theorem (see e.g. Zeidler (1985), p.130) and the first-order condition $\widetilde{\mathcal{L}}^{(1)}({\bf B}^*(\theta);\theta)=0$ yields that ${\bf B}^{*}( \theta)$ is continuously differentiable in $\theta$ over  $\widetilde{\Theta}^0$, with
\begin{equation}
\frac{
\partial
{\bf B}^{*}(\theta)
}
{
\partial \theta^T
}
=
-
\left[
\widetilde{\mathcal{L}}^{(2)} ({\bf B}^* (\theta);\theta)
\right]^{-1}
\frac{\partial
\widetilde{\mathcal{L}}^{(1)} ({\bf B}^* (\theta);\theta)
}{\partial \theta^T}
.
\label{Bstarder}
\end{equation}
Recall now that $\Theta^0 \subset \widetilde{\Theta}^0$ when $\overline{h}$ tends to $0$. Hence continuity of ${\bf B}^{*}(\cdot)$, $\partial \widetilde{\mathcal{L}}^{(1)} (\cdot,\cdot)/\partial \theta^{T}$ and compactness of $\widetilde{\Theta}^0$ give
\begin{eqnarray}
\lim_{\overline{h} \rightarrow 0}
\sup_{\theta \in \Theta^0}
  \left\|
    {\bf B}^{*}(\theta)
      -
    {\bf B}^{*}(\alpha;0,x)
  \right\|
& = &
0
,
\nonumber
\\
\lim_{ \overline{h} \rightarrow 0} \sup_{\theta \in \Theta^0}
\left\|
  \frac{\partial \widetilde{\mathcal{L}}^{(1)}({\bf B}^{*}(\theta);\theta)}{\partial \theta^T}
    -
  \frac{\partial \widetilde{\mathcal{L}}^{(1)}({\bf B}^{*}(\alpha;0,x);\alpha,0,x)}{\partial \theta^T}
\right\|
& = &
0 .
\label{DiffGtilde1}
\end{eqnarray}
Since the first limit is (\ref{B*2B*0}), (ii) is proved.

{\it Proof of (iii)}.
We bound the partial derivative (\ref{Bstarder}).
Observe that (\ref{B*2B*0}), the expression of ${\bf B}^{*}(\alpha;0,x)$, the compactness of $\Theta^0$ and Assumption \ref{F} yield that there is a compact $\mathcal{B}$ such that
${\bf B}^* (\theta)$ is in $\mathcal{B}$ for all $\theta$ in $\Theta^0$, provided $\overline{h}$ is small enough. Then (\ref{Gtilde2}) and (\ref{Formquad}) give that  uniformly in $\theta$ in $\Theta^0$,
\[
\widetilde{\mathcal{L}}^{(2)}({\bf B}^{*}(\theta);\theta) \succ C \int_{\mathcal{B}(0,1)} {\bf U}(z){\bf U}(z)^T dz
.
\]
Hence
(\ref{Bstarder}) and (\ref{DiffGtilde1}) give
\begin{equation}
\lim_{\overline{h} \rightarrow 0}
\sup_{\theta \in \Theta^0}
\left\|
  \frac{
    \partial
       {\bf B}^{*}(\theta)
       }
       {
       \partial \theta^T
       }
\right\|
\leq
C
\left\|
    \left(
       \int_{\mathcal{B}(0,1)}
         {\bf U}(z)
         {\bf U}(z)^T
       dz
    \right)^{-1}
\right\|
\lim_{\overline{h} \rightarrow 0}
\sup_{\theta \in \Theta^0}
\left\|
  \frac{
    \partial
     \widetilde{\mathcal{L}}^{(1)}
       ({\bf B}^{*}(\theta);\theta)
       }
       {
       \partial \theta^T
       }
\right\|
\leq
C
.
\label{BornenormBetoile}
\end{equation}

Let us now return to the proof of (iii). The differentiability results above yield that
$\theta \in \Theta^1 \mapsto Q^{*}(x';\theta) = {\bf U}((x -x')/h)^T {\bf B}^{*}(\theta)$ is continuously differentiable in $\theta$. We have for all $x$, $x'$ in $\mathcal{X}$ and $h \geq \underline{h}$,
\[
\left\|
{\bf U}\left(
\frac{x-x'}{h}
\right)
\right\|
\leq
\frac{C}{\underline{h}^p}
,
\quad
\left\|
\frac{\partial}{\partial \theta^T}
{\bf U}\left(
\frac{x-x'}{h}
\right)
\right\|
\leq
\frac{C}{\underline{h}^{p+1}}
.
\]
Hence for $\overline{h}$ small enough, (\ref{B*2B*0}) and (\ref{BornenormBetoile}) yield that for all $\theta$ in $\Theta^1$ and $x'$ in $\mathcal{X}$,
\begin{eqnarray*}
\left\|
\frac{\partial Q^{*}(x';\theta) }{\partial \theta^T}
\right\|
& = &
\left\|
\left[
\frac{\partial}{\partial \theta^T}{\bf U}\left(\frac{x-x'}{h}\right)^{T}
\right]
{\bf B}^{*}(\theta)
+
{\bf U}\left(\frac{x-x'}{h}\right)^T
\frac{\partial {\bf B}^{*}(\theta) }{\partial \theta^T}
\right\|
\\
& \leq &
\left\|
\frac{\partial}{\partial \theta^T}{\bf U}\left(\frac{x-x'}{h}\right)
\right\|
\left\|
{\bf B}^{*}(\theta)
\right\|
+
\left\|
{\bf U}\left(\frac{x-x'}{h}\right)
\right\|
\left\|
\frac{\partial {\bf B}^{*}(\theta)}{\partial \theta^T}
\right\|
\leq
C\underline{h}^{-p}\left( 1 + \underline{h}^{-1}\right)
.
\end{eqnarray*}
The Taylor inequality shows that (iii) is proved.

{\it Proof of (iv).} The change of variable $x' = x+hz$ shows that it is sufficient to prove that, for all $\theta$ in $\Theta^0$ and $z$ in $\mathcal{K}$,
\[
f(Q^{*}(x + hz;\theta) |x +hz)
\geq
C
\mbox{ \rm with }
f(Q^{*}(x+hz;\theta)|x+hz) = f({\bf U}(z)^T {\bf B}^{*}(\theta)| x+hz)
,
\]
which is true for $\overline{h}$ small enough by (\ref{B*2B*0}) and under Assumption \ref{F} which gives that $f(y|x) \geq C > 0$ for $y$ in any compact subset of $\Rit$ and any $x$ in $\mathcal{X}_0$.
\eop

\subsection{Proof of Proposition \protect{\ref{R1order}}}
The proof of the Proposition uses the two following Lemmas. In what follows, the stochastic processes $R (\cdot;\cdot)$, ${\bf R}^1 (\cdot;\cdot)$ and ${\bf R}^2 (\cdot;\cdot)$ have the same distribution than the $R_i (\cdot;\cdot)$, ${\bf R}^1 (\cdot;\cdot)$ and ${\bf R}^2 (\cdot;\cdot)$ in (\ref{Rho}), (\ref{R1}) and (\ref{R2}).
Define also
\begin{equation}
\delta (\beta,\theta) 
=
{\bf U}\left(\frac{X-x}{h}\right)^{T}
\frac{\beta}{\left(n h^{d}\right)^{1/2}}
.
\label{Delta}
\end{equation}
\begin{lem}
Under Assumptions \ref{F}, \ref{K} and \ref{X}, we have 
$$
\Var
\left(
R (\beta,\epsilon;\theta)
\right)
\leq
C
\frac{\left\|\epsilon\right\|^2(\left\|\beta\right\|+ \left\|\epsilon\right\|)}
     {
  n
  \left(n h^{d}\right)^{1/2}
     }
.
$$
\label{Bornesrho1}
\end{lem}
\noindent
{\bf Proof of Lemma \ref{Bornesrho1}.}
Observe  $\ell_{\alpha} (t)= 2 \int_{0}^{t} (\alpha - \ind (z \leq 0))dz$. Hence (\ref{Rho}) and (\ref{Delta}) yield 
\begin{equation}
R
\left(\beta,\epsilon;\theta\right)
= 
2K_h(X-x)
\int_{\delta (\beta,\theta)}^{\delta (\beta,\theta) + \delta (\epsilon,\theta)}
\left(
\ind\left(Y\leq Q^{*}(X;\theta) + t\right)
-
\ind\left(Y\leq Q^{*}(X;\theta) \right)
\right)
dt.
\label{Rho0.2}
\end{equation}
The Cauchy-Schwarz inequality 
give
\begin{eqnarray*}
R
\left(\beta,\epsilon;\theta\right)^2
& = & 
4K_h(X-x)^2
\left(
  \int_{\delta (\beta,\theta)}^{\delta (\beta,\theta) +     
  \delta (\epsilon,\theta)}
    \left(
\ind\left(Y\leq Q^{*}(X;\theta) + t\right)
-
\ind\left(Y\leq Q^{*}(X;\theta) \right)
\right)
   dt
\right)^2
\\
& \leq &
4K_h(X-x)^2
 \left|\delta (\epsilon,\theta) \right|
\left|
  \int_{\delta (\beta,\theta)}^{
    \delta (\beta,\theta) + \delta (\epsilon,\theta)}
\left(
\ind\left(Y\leq Q^{*}(X;\theta) + t\right)
-
\ind\left(Y\leq Q^{*}(X;\theta) \right)
\right)^2
dt
\right|
\\
& < &
4
K_h(X-x)^2
\left| \delta (\epsilon,\theta) \right|
\left|
\int_{\delta (\beta,\theta)}^{
\delta (\beta,\theta) + \delta (\epsilon,\theta)}
\ind\left(
\left|
Y
-
Q^{*}(X;\theta)
\right| 
<
\left|t \right|
\right)
dt
\right|
.
\end{eqnarray*}
Hence Assumption \ref{F} and (\ref{Delta}) give
\begin{eqnarray*}
\esp
\left[
R^{2}
\left(\beta,\epsilon;\theta\right)
\left| X \right.
\right]
& \leq  &
4
K_{h}\left(X-x\right)^{2}
\left|
\delta (\epsilon,\theta)
\right|
\left|
\int_{\delta (\beta,\theta)}^{\delta (\beta,\theta) + \delta (\epsilon,\theta)}
\left\{
\int
\ind\left(
\left|
y-
Q^{*}(X;\theta)
\right| 
<
\left|t \right|
\right)
f(y|X)
dy
\right\}
dt
\right|
\\
& \leq &
4
K_{h}\left(X-x\right)^{2}
\left\|
f(\cdot|\cdot)
\right\|_{\infty}
|
  \delta (\epsilon,\theta)
|
\left|
2
\int_{\delta (\beta,\theta)}^{\delta (\beta,\theta) + \delta (\epsilon,\theta)}
|t|
dt
\right|
\\
& \leq &
C
K_{h}\left(X-x\right)^{2}
\delta (\epsilon,\theta)^2
  \left(|\delta (\beta,\theta)| +  
  |\delta (\epsilon,\theta)|
\right)
\\
& \leq &
C
\frac{K^2
\left(\frac{X-x}{h}\right)
\left\|
{\bf U}
\left(\frac{X-x}{h}\right)
\right\|^3
}{
\left(nh^d\right)^{3/2}}
\left\|\epsilon\right\|^2
\left(\left\|\beta\right\|+\left\|\epsilon\right\|\right)
.
\end{eqnarray*}
Then, under Assumptions \ref{K} and \ref{X},
\begin{eqnarray*}
\Var
\left(R\left(\beta,\epsilon;\theta\right)\right)
& \leq & 
\esp [R^2 \left(\beta,\epsilon;\theta\right)]
= 
\esp
\left[\esp [R^2 \left(\beta,\epsilon;\theta\right)|X]\right] 
\\
& \leq &
\frac{C
\left\|\epsilon\right\|^2
\left(\left\|\beta\right\|+\left\|\epsilon\right\|\right)
}
{\left(nh^{d}\right)^{3/2}}
\int
K^2
\left(\frac{x'-x}{h}\right)
\left\|
{\bf U}
\left(\frac{x'-x}{h}\right)
\right\|^3f_{X}\left(x' \right)
d x'
\\
& 
\leq 
&
\frac{C
\left\|\epsilon\right\|^2
\left(\left\|\beta\right\|+\left\|\epsilon\right\|\right)
}
{\left(nh^{d}\right)^{3/2}}
h^d
\int
K^2
\left(z\right)
\left\|
{\bf U}
\left(z\right)
\right\|^3f_{X}\left(x+hz \right)
dx'
\leq 
C
\frac{
\left\|\epsilon\right\|^{2}
\left(\left\|\beta\right\| + \left\|\epsilon\right\| \right)
}
{n \left( nh^d \right)^{1/2}}
.
\eop
\end{eqnarray*}

Define
\[
\mathcal{F}
=
\mathcal{F}
\left(t_{\beta},t_{\epsilon},\Theta^1 \right)
=
\left\{
  R
    \left(
      \beta , \epsilon; \theta
    \right)
    ,
      \left(\beta,\epsilon,\theta\right)
        \in
      \mathcal{B}(0,t_{\beta}) \times \mathcal{B}(0,t_{\epsilon})
      \times
      \Theta^1
\right\}
 .
\]
The next lemma studies coverings of $\mathcal{F}$ with brackets $\left[\underline{R},\overline{R} \right]$. Recall that the bracket
$\left[\underline{R},\overline{R} \right]= \left[\underline{R} (X,Y),\overline{R}(X,Y) \right]$ is the set of random variables
$r=r(X,Y)$ such that
$\underline{R} \leq r \leq \overline{R}$ almost surely. 
\begin{lem} 
Under Assumptions \ref{F}, \ref{K} and \ref{X} and if $t_{\beta}+t_{\epsilon} \geq 1$ and $n$ is large enough,
\begin{enumerate}
\item
There are some  $\overline{\sigma}^2$ and $\overline{w}$, with
$$
\overline{\sigma}^2
\asymp
\frac{t_{\epsilon}^2
(
 t_{\epsilon}+t_{\beta}
)}{n(n\underline{h}^d)^{1/2}}
,
\quad
\overline{w}
\asymp
\frac{t_{\beta} + t_{\epsilon}}{(n\underline{h}^d)^{1/2}}
,
$$
such that for all integer number $k  \geq  2$, $(\beta,\epsilon,\theta)$ 
in 
$\mathcal{B}(0,t_{\beta}) 
\times 
\mathcal{B}( 0,t_{\epsilon})
\times
\Theta^1$,
\begin{eqnarray*}
\esp
\left[
\left|
R (\beta,\epsilon;\theta)
-
\esp
\left[
R (\beta,\epsilon;\theta)
\right]
\right|^{k}
\right]
\leq
\frac{k!}{2}
\overline{w}^{k-2}
\overline{\sigma}^2
.
\end{eqnarray*}
\item
Let $\tau$ in $(0,1)$ be a bracket length. There is an  set of brackets 
$\mathcal{I}_{\tau}=\left\{\left[\underline{R}_{j,\tau},\overline{R}_{j,\tau}\right], 1 \leq j \leq e^{H(\tau)} \right\}$ such that
\begin{eqnarray*}
&&
\mathcal{F} 
\subset
\bigcup_{1 \leq j \leq e^{H(\tau)}}
\left[\underline{R}_{j,\tau},\overline{R}_{j,\tau}\right],
\\
&&
\esp\left[
  \left|
    \underline{R}_{j,\tau} 
      -
    \overline{R}_{j,\tau} 
  \right|^k
  \right]
\leq 
\frac{k!}{2}
\overline{w}^{k-2}
\tau^2
\mbox{ \it for all integer number $k \geq 2$ and all $j$ in $\left[1,e^{H(\tau)}\right]$,}
\\
&&
H (\tau)  \leq 
C \log \left( 
\frac{n(t_{\beta} + t_{\epsilon}) }{\tau}\right)
\mbox{ \it for all $\tau$, $t_{\beta}$ and $t_{\epsilon}$.}  
\end{eqnarray*}
\end{enumerate}
\label{Controlerho1}
\end{lem}
\noindent
{\bf Proof of Lemma \ref{Controlerho1}}.
Define for $\beta$ in $\Rit^P$ 
$$
\widetilde{R} (\beta;\theta) 
= 
2K_h (X-x) \int_0^{\delta (\beta,\theta)} 
\left( 
  \ind \left(Y \leq Q^{*}(X;\theta) +u\right)
    -
  \ind \left(Y \leq Q^{*}(X;\theta)\right)
\right)
du
.
$$
Let $\mathrm{sgn}(t) = \ind (t \geq 0) - \ind (t<0)$. Observe that
$\widetilde{R} (\beta;\theta) \geq 0$ with
\begin{eqnarray}
\widetilde{R} (\beta;\theta) 
& = & 
2K_h (X-x) 
\int_0^{|\delta (\beta,\theta)|} 
\left|
  \ind \left(Y \leq Q^{*}(X;\theta) +\mathrm{sgn}(\delta (\beta,\theta))u\right)
    -
  \ind \left(Y \leq Q^{*}(X;\theta)\right)
\right|
du
\nonumber \\
& = & 
2K_h (X-x) 
|\delta (\beta,\theta)|
\int_0^{1} 
\left|
  \ind \left(Y \leq Q^{*}(X;\theta) +\delta (\beta,\theta)v\right)
    -
  \ind \left(Y \leq Q^{*}(X;\theta)\right)
\right|
dv
\nonumber \\
& = &
2K_h (X-x) 
|\delta (\beta,\theta)|
\int_0^{1} 
\ind
\left(
\mbox{\rm
$Y-Q^{*}(X;\theta)$ lies between 
$0$ and $\delta (\beta,\theta)v$}
\right)
dv
.
\label{TildeR}
\end{eqnarray}
(\ref{Rho0.2}) and $\delta (\beta,\theta) + \delta (\epsilon,\theta) = \delta (\beta+\epsilon,\theta)$ give
\begin{eqnarray}
R(\beta;\epsilon,\theta)
& = &
\widetilde{R}(\beta + \epsilon;\theta)
-
\widetilde{R}(\beta;\theta)
\label{Decomposition}
.
\end{eqnarray}
It also follows from (\ref{Delta}) and Assumption \ref{K} that for all $\beta$ in $\mathcal{B} \left(0,t_{\beta}+t_{\epsilon} \right)$ and all $\theta$ in $\Theta^1$
\begin{equation}
\left|
\widetilde{R}(\beta ;\theta)
\right|
\leq
2
\left\|{\bf U} \left(\frac{X-x}{h}\right)\right\|
K\left(\frac{X-x}{h}\right)
\frac{\left\|\beta\right\|}{\left(n h^d\right)^{1/2}}
\leq
\frac{\overline{w}}{2},
\quad
\overline{w}
\asymp
\frac{t_{\beta}+t_{\epsilon}}{\left(n \underline{h}^d\right)^{1/2}}.
\label{BdtildeR}
\end{equation}
Part (i) follows from Lemma \ref{Bornesrho1} and (\ref{Decomposition}) which give
\begin{eqnarray*}
\lefteqn{
\esp
\left[
\left|
R (\beta,\epsilon;\theta)
-
\esp
\left[
R (\beta,\epsilon;\theta)
\right]
\right|^{k}
\right]
}
&&
\\
& = &
\esp
\left[
\left|
\widetilde{R} (\beta+\epsilon;\theta)
-
\esp
\left[
\widetilde{R} (\beta+\epsilon;\theta)
\right]
-
\left(
\widetilde{R} (\beta;\theta)
-
\esp
\left[
\widetilde{R} (\beta;\theta)
\right]
\right)
\right|^{k-2}
\left|
R (\beta,\epsilon;\theta)
-
\esp
\left[
R (\beta,\epsilon;\theta)
\right]
\right|^{2}
\right]
\\
& \leq &
\left(2 \times \frac{\overline{w}}{2} \right)^{k-2}
\Var\left(R (\beta,\epsilon;\theta)\right)
\leq
\overline{w}^{k-2}
\overline{\sigma}^2
.
\end{eqnarray*}
The proof of part (ii) will be divided in three steps. 
Let 
$
\widetilde{\mathcal{F}}_{t}
$
be
$ 
\{
\widetilde{R}(\beta;\theta)
,
(\beta,\theta)
\in
\mathcal{B}(0,t)
\times
\Theta^1
\}
$
.
For the sake of brevity we abbreviate $\underline{R}_{j,\tau},\overline{R}_{j,\tau}$ into $\underline{R}_{j},\overline{R}_{j}$.

{\it Step 1 : Coverings of $\mathcal{F}$ and $\widetilde{\mathcal{F}}_{t}$,  $t=t_{\beta}+t_{\epsilon} \geq 1$.}   We show in this step that it is sufficient to find a covering of $\widetilde{\mathcal{F}}_{t}$ with 
$H(\tau)=H(\tau;t)$ brackets satisfying
\begin{eqnarray}
\esp
\left[
\left|
\underline{R}_{j}
  -
\overline{R}_{j}
\right|^k
\right]
& \leq &
\frac{k!}{8}
\left(\frac{\overline{w}}{2}\right)^{k-2}
\tau^2
,
\label{TBPbracks}
\\
H(t)
& \leq &
C \log\left(\frac{nt}{\tau}\right)
.
\label{TBPH}
\end{eqnarray}
Indeed, consider two such coverings of $\widetilde{\mathcal{F}}_{t_{\beta}}$ and  
$\widetilde{\mathcal{F}}_{t_{\beta}+t_{\epsilon}}$,
$$
\widetilde{\mathcal{F}}_{t_{\beta}}
\subset
\bigcup_{1 \leq j \leq e^{H_1(\tau)}}
\left[\underline{R}_{j}^{1},\overline{R}_{j}^{1} \right]
,
\quad
\widetilde{\mathcal{F}}_{t_{\beta}+t_{\epsilon}}
\subset
\bigcup_{1 \leq j \leq e^{H_2(\tau)}}
\left[\underline{R}_{j}^{2},\overline{R}_{j}^{2} \right],
$$
$H_1 ( \tau) \leq H_2 ( \tau) = H (\tau, t_{\beta}+t_{\epsilon})$.
Consider a $R(\beta,\epsilon;\theta)$ in $\mathcal{F}$. Since 
$
\widetilde{R}(\beta;\theta)
\in
\left[\underline{R}_{j_1}^{1},\overline{R}_{j_1}^{1} \right]
$
and
$
\widetilde{R}(\beta+\epsilon;\theta)
\in
\left[\underline{R}_{j_2}^{2},\overline{R}_{j_2}^{2} \right]
$ for some $j_1$ and $j_2$, (\ref{Decomposition}) implies that
$
R (\beta,\epsilon;\theta)
\in
\left[
\underline{R}_{j_2}^{2}-\overline{R}_{j_1}^{1},\overline{R}_{j_2}^{2}-\underline{R}_{j_1}^{1} \right]
$.
Hence these $e^{H'(\tau)}$ brackets form a covering of $\mathcal{F}$
with, using (\ref{TBPbracks}) and (\ref{TBPH}),
\begin{eqnarray*}
\esp
\left[
\left|
\overline{R}_{j_2}^{2}-\underline{R}_{j_1}^{1}
-
\left(\overline{R}_{j_2}^{2}-\underline{R}_{j_1}^{1}\right)
\right|^k
\right]
& \leq &
2^{k-1}
\left(
\esp
\left[
\left|
\overline{R}_{j_2}^{2}
-
\underline{R}_{j_2}^{2}
|\right|^k
\right]
+
\esp
\left[
\left|
\overline{R}_{j_1}^{1}
-
\underline{R}_{j_1}^{1}
|\right|^k
\right]
\right)
\\
&\leq &
2^k
\frac{k!}{8}
\left(\frac{\overline{w}}{2}\right)^{k-2}
\tau^2
=
\frac{k!}{2}
\overline{w}^{k-2}
\tau^2
,
\\
H'(\tau)=H_1 (\tau)+H_2 (\tau)
& \leq &
C \log\left(\frac{n\left(t_{\beta}+t_{\epsilon}\right)}{\tau}\right)
.
\end{eqnarray*}

{\it Step 2: Preliminary results for the construction of a covering of  $\widetilde{\mathcal{F}}_{t}$.} 
We bound the increments of $(\beta,\theta) \mapsto Q^*(X;\theta), K_h (X-x), \delta(\beta,\theta)$.
Lemma \ref{Biais1}-(iii) gives that for all $\theta$, $\theta'$ in $\Theta^1$
$$
\left|
Q^{*} \left(X ; \theta \right) - Q^{*} \left(X ; \theta' \right)
\right|
\leq
C
\underline{h}^{-p}(1+\underline{h}^{-1})
\left\|
\theta - \theta'
\right\|
.
$$
Under Assumption \ref{K}
\begin{eqnarray*}
\left|K\left(\frac{X-x}{h}\right) - K\left(\frac{X-x'}{h'}\right)\right|
& \leq &
C
\left(
    \left\|
      \frac{x - x'}{h'}
    \right\|
+
    \left\|
      X-x
    \right\|
    \left|
      \frac{1}{h}
        -
      \frac{1}{h'}
    \right|
\right)
\\
& \leq &
  C \left(
    \frac{1}{\underline{h}} \left\|x - x'\right\|
      +
    \frac{1}{\underline{h}^2} \left|h - h' \right|
    \right)
\leq
\frac{C}{\underline{h}^2}
\left\|
  \theta - \theta'
\right\|
\nonumber
.
\end{eqnarray*}
For the increments of $\delta (\beta,\theta)$, define
${\bf U} = {\bf U} ( X - x )$, ${\bf U}' = {\bf U} ( X - x' )$,  ${\bf H}' =
{\bf H}(h')$. This gives
\begin{eqnarray*}
\lefteqn{
\left|
  \delta (\beta,\theta) - \delta(\beta',\theta')
\right|
}
&&
\\
& = &
\left|
   {\bf U}^T \frac{{\bf H}^{-1}}{\left(nh^d\right)^{1/2}} \left( \beta - \beta' \right)
   +
   \left({\bf U}' - {\bf U}\right)^T \frac{{\bf H}^{-1}}{\left(nh^d\right)^{1/2}} \beta'
   +
   {\bf U}^{\prime T} \left(\frac{{\bf H}^{-1}}{\left(nh^d\right)^{1/2}} -
\frac{{\bf H}^{\prime -1}}{\left(nh^{\prime d}\right)^{1/2}}\right) \beta'
   \right|
\\
& \leq &
C
\left\|
  \frac{{\bf H}^{-1}}{\left(nh^d\right)^{1/2}}
\right\|
\left\|
  \beta - \beta'
\right\|
+
\left\| x-x'\right\|
\left\|
  \frac{{\bf H}^{-1}}{\left(nh^d\right)^{1/2}}
\right\|
\left\|
  \beta'
\right\|
+
C
\left\| \beta' \right\|
\left\|
  \frac{{\bf H}^{-1}}{\left(nh^d\right)^{1/2}}
    -
  \frac{{\bf H}^{\prime-1}}{\left(nh^{\prime d}\right)^{1/2}}
\right\|
\\
& \leq &
\frac{C(1+t)}{\underline{h}^{p} \left(n\underline{h}^d\right)^{1/2}}
\left(
  \left\|
    \beta - \beta'
  \right\|
+
  \left\|
    x - x'
  \right\|
+
\frac{1}{
  \underline{h}}
  \left|
    h - h'
  \right|
\right)
\leq
\frac{C(1+t)}{\underline{h}^{p+1} \left(n\underline{h}^d\right)^{1/2}}
\left(
\left\|
    \beta - \beta'
  \right\|
+
\left\|
  \theta - \theta'
\right\|
\right)
.
\end{eqnarray*}

{\it Step 3 : Construction of the covering of $\widetilde{\mathcal{F}}_{t}$.}
Define
\begin{eqnarray*}
\rho (q,\delta)
& = &
\left|
\ind \left(q \leq \delta \right)
-
\ind \left(q \leq 0 \right)
\right|
=
\ind
\left(
q \in (0,\delta ]
\right)
\ind
\left(
\delta \geq 0
\right)
+
\ind
\left(
q \in [\delta,0)
\right)
\ind
\left(
\delta < 0
\right)
,
\\
r(q,\delta)
&=&
\int_0^1
\rho (q,\delta v)dv
.
\end{eqnarray*}
Hence (\ref{TildeR}) shows
$$
\widetilde{R} (\beta;\theta)
=
2K_h (X-x) 
|\delta (\beta,\theta)|
r
\left(
Y-Q^* (X;\theta)
,
\delta
\left(
\beta
,
\theta
\right)
\right)
.
$$
For any $\eta>0$, there exists functions $\underline{\rho} (q,\delta)=\underline{\rho}_{\eta} (q,\delta)$ and 
$\overline{\rho} (q,\delta)=\overline{\rho}_{\eta} (q,\delta)$ and an open set $D=D_{\eta} \subset \Rit^2$ such that
$$
\begin{array}{ll}
\rho-{\rm (i)}
&
0 
\leq 
\underline{\rho} (q,\delta) 
\leq 
\rho (q,\eta) 
\leq 
\overline{\rho} (q,\delta)
\leq 1
\mbox{ \rm
for all $(q,\delta)$,
with
}
\underline{\rho} (q,\delta) 
= 
\rho (q,\eta) 
= 
\overline{\rho} (q,\delta)
\mbox{ \rm
if
} 
(q,\delta) \in \Rit^2\setminus D_{\eta}
,
\\
\rho-{\rm (ii)}
&
\sup_{(q,\delta) \in D_{\eta}}
\left(
\left|
\frac{\partial \underline{\rho}  (q,\delta)}{\partial q}
\right|
+
\left|
\frac{\partial \underline{\rho}  (q,\delta)}{\partial \delta}
\right|
+
\left|
\frac{\partial \overline{\rho}  (q,\delta)}{\partial q}
\right|
+
\left|
\frac{\partial  \overline{\rho}  (q,\delta)}{\partial \delta}
\right|
\right)
\leq
C \eta^{-1/2},
\\
\rho-{\rm (iii)}
&
D
\subset
D'=
\left\{
\left(q,\delta\right)
\in 
\Rit^2
;
|q| \leq C \eta^{-1/2}
\mbox{ \rm or }
\left|q-\delta \right|
\leq C \eta^{-1/2}
\right\}
.
\end{array}
$$
Define
$\underline{r} (q,\delta) = \int_0^1 \underline{\rho} (q,v\delta) dv$,
$\overline{r} (q,\delta) = \int_0^1 \overline{\rho} (q,v\delta) dv$ and
\begin{eqnarray*}
\underline{R} (\beta,\theta)
& = &
2K_h (X-x) 
|\delta (\beta,\theta)|
\underline{r}
\left(
Y-Q^* (X;\theta)
,
\delta
\left(
\beta
,
\theta
\right)
\right)
,
\\
\overline{R} (\beta,\theta)
& = &
2K_h (X-x) 
|\delta (\beta,\theta)|
\overline{r}
\left(
Y-Q^* (X;\theta)
,
\delta
\left(
\beta
,
\theta
\right)
\right)
.
\end{eqnarray*}
Since $K(\cdot) \geq 0$, $\rho$-(i) gives that these functions are such that
\begin{equation}
\underline{R} (\beta,\theta)
\leq
\widetilde{R} (\beta,\theta)
\leq
\overline{R} (\beta,\theta)
.
\label{Brack1}
\end{equation}
We now bound $\underline{R} (\beta,\theta) - \underline{R} (\beta',\theta')$ and 
$\overline{R} (\beta,\theta) - \overline{R} (\beta',\theta')$. We have
\begin{eqnarray*}
\lefteqn{
\left|
\underline{R} (\beta,\theta) - \underline{R} (\beta',\theta')
\right|
\leq 
2
\left|
K_h (X-x)
-
K_{h'} (X-x')
\right|
|\delta (\beta,\theta)|
\underline{r}
\left(
Y-Q^* (X;\theta)
,
\delta
\left(
\beta
,
\theta
\right)
\right)
}
\\
&&
+ 
2
K_{h'} (X-x')
|\delta (\beta,\theta) - \delta (\beta',\theta')|
\underline{r}
\left(
Y-Q^* (X;\theta)
,
\delta
\left(
\beta
,
\theta
\right)
\right)
\\
&&
+ 
2
K_{h'} (X-x')
| \delta (\beta',\theta')|
\left|
\underline{r}
\left(
Y-Q^* (X;\theta)
,
\delta
\left(
\beta
,
\theta
\right)
\right)
-
\underline{r}
\left(
Y-Q^* (X;\theta')
,
\delta
\left(
\beta'
,
\theta'
\right)
\right)
\right|.
\end{eqnarray*}
Hence Step 1, $\rho$-(i,ii), (\ref{Delta}) and the Taylor inequality give for all $(\beta,\theta)$,  $(\beta',\theta')$ in $\mathcal{B} (0,t) \times \Theta^1$, provided $n$ is large enough,
\begin{eqnarray*}
\left|
\underline{R} (\beta,\theta) - \underline{R} (\beta',\theta')
\right|
&\leq&
C
\left[
\frac{t}{\underline{h}^p (n\underline{h}^d)^{1/2}}
\frac{\left\|\theta-\theta'\right\|}{\underline{h}^2}
+
\frac{1+t}{\underline{h}^{p+1} (n\underline{h}^d)^{1/2}}
\left(
\left\|\theta-\theta'\right\|
+
\left\|\beta-\beta'\right\|
\right)
\right]
\\
&&
+
C\eta^{-1/2}
\left[
\frac{\left\|\theta-\theta'\right\|}{\underline{h}^{p+1}}
+
\frac{1+t}{\underline{h}^{p+1} (n\underline{h}^d)^{1/2}}
\left(
\left\|\theta-\theta'\right\|
+
\left\|\beta-\beta'\right\|
\right)
\right]
\\
& \leq &
C
\frac{\left(1+\eta^{-1/2}\right)\left(1+t\right)}{\underline{h}^{p+2} }
\left(
\left\|\theta-\theta'\right\|
+
\left\|\beta-\beta'\right\|
\right).
\end{eqnarray*}
Arguing symmetrically gives
$$
\left|
\overline{R} (\beta,\theta) - \overline{R} (\beta',\theta')
\right|
\leq
C
\frac{\left(1+\eta^{-1/2}\right)\left(1+t\right)}{\underline{h}^{p+2} }
\left(
\left\|\theta-\theta'\right\|
+
\left\|\beta-\beta'\right\|
\right).
$$

We now construct the brackets. Recall that there is a covering of $\mathcal{B}(0,t) \times \Theta^1$ with $N$ balls $\mathcal{B}\left((\beta_j,\theta_j),\eta\right)$, $\theta_j = (\alpha_j,h_j,x_j)$, with center $(\beta_j,\theta_j)$ and radius $\eta$ such that
\begin{equation}
N \leq
\max
\left(
1,
\frac{Ct^P}{\eta^{P+d+2}}
\right)
,
\label{Neta}
\end{equation}
see van de Geer (1999, p.20).
Define
$$
\underline{R}_j'
=
\underline{R} (\beta_j,\theta_j)
-
C
\eta
\frac{\left(1+\eta^{-1/2}\right)\left(1+t\right)}{\underline{h}^{p+2} }
,
\quad
\overline{R}_j'
=
\overline{R} (\beta_j,\theta_j)
+
C
\eta
\frac{\left(1+\eta^{-1/2}\right)\left(1+t\right)}{\underline{h}^{p+2} },
$$
\begin{equation}
\underline{R}_j=\max \left(0,\underline{R}_j'\right),
\quad
\overline{R}_j=\min\left(\frac{\overline{w}}{2},\overline{R}_j'\right)
.
\label{Brackets}
\end{equation}
Bounding $\underline{R} (\beta,\theta)-\underline{R}_j$ and $\overline{R} (\beta,\theta)-\overline{R}_j$ for $(\beta,\theta)$ in $\mathcal{B}\left((\beta_j,\theta_j),\eta\right)$, (\ref{Brack1}) and (\ref{BdtildeR}) give   
\begin{equation}
\underline{R}_j'
\leq
\underline{R}_j
\leq
\widetilde{R} (\beta,\theta)
\leq
\overline{R}_j
\leq
\overline{R}_j'.
\label{Brack2}
\end{equation}
It then follows that
$\left\{\left[\underline{R}_j,\overline{R}_j\right], j=1,\ldots,N\right\}$
is a covering of $\widetilde{\mathcal{F}}_t$ with, since $0 \leq \underline{R}_j \leq \overline{R}_j \leq \overline{w}/2$,
\begin{equation}
\left|
  \overline{R}_{j}
-
  \underline{R}_{j}
\right|
 \leq 
\frac{\overline{w}}{2} \asymp
C\frac{t}{ \left(n\underline{h}^d\right)^{1/2}}
.
\label{Rhobracket}
\end{equation}

We now bound 
$\esp 
\left[
\left(
\overline{R}_{j}
-
\underline{R}_{j}
\right)^2
\right]$
and
$\esp 
\left[
\left|
\overline{R}_{j}
-
\underline{R}_{j}
\right|^k
\right]$.
(\ref{Brack2}), $\rho$-(i,iii), (\ref{Delta}) and Assumptions \ref{F}, \ref{K} give
\begin{eqnarray*}
\lefteqn{
\esp 
\left[
\left(
\overline{R}_{j}
-
\underline{R}_{j}
\right)^2
\right]
\leq
\esp 
\left[
\left(
\overline{R}_{j}'
-
\underline{R}_{j}'
\right)^2
\right]
\leq
2
\esp 
\left[
\left(
\overline{R}\left(\beta_j,\theta_j\right)
-
\underline{R}\left(\beta_j,\theta_j\right)
\right)^2
\right]
+
C
\eta^2
\frac{\left(1+\eta^{-1/2}\right)^2\left(1+t\right)^2}{\underline{h}^{2(p+2)} }
}
&&
\\
& \leq &
8
\esp
\left[
K_{h_j}^2\left(X-x_j\right)
\delta^2\left(\beta_j,\theta_j\right)
\left(
\overline{r}
\left(
Y-Q^*\left(X;\theta_j\right)
,
\delta\left(\beta_j,\theta_j\right)
\right)
-
\underline{r}
\left(
Y-Q^*\left(X;\theta_j\right)
,
\delta\left(\beta_j,\theta_j\right)
\right)
\right)^2
\right]
\\
&&
+
C
\frac{\left(1+t\right)^2}{ \underline{h}^{2(p+2)}}
\left(\eta^2+\eta\right)
\\
& \leq &
8
\esp
\left[
K_{h_j}^2\left(X-x_j\right)
\delta^2\left(\beta_j,\theta_j\right)
\int
\left(
\int_0^1
\ind
\left(
\left(y-Q^*(X;\theta_j),v \delta\left(\beta_j,\theta_j\right) \right)
\in D
\right)
dv
\right)^2
f(y|X)
dy
\right]
\\
&&
+
C
\frac{\left(1+t\right)^2}{ \underline{h}^{2(p+2)}}
\left(\eta^2+\eta\right)
\\
& \leq &
\frac{8 h_j^d \left\|\beta\right\|^2}{n h_j^d}
\int
K^2 (z)
\left\|\mathbf{U} (z) \right\|^2
\\
&&
\times
\left[
\int
\int_0^1
\ind
\left(
\left(y-Q^*(x_j+h_jz;\theta_j),v \delta\left(\beta_j,\theta_j\right) \right)
\in D
\right)
dv
f(y|x_j+h_j z)
dy
\right]
f(x_j+h_j z) dz
\\
&&
+
C
\frac{\left(1+t\right)^2}{ \underline{h}^{2(p+2)}}
\left(\eta^2+\eta\right)
\\
& \leq &
C
\frac{\left(1+t\right)^2}{ \underline{h}^{2(p+2)}}
\left(\eta^2+\eta+\eta^{1/2} \right).
\end{eqnarray*}
This together with (\ref{Rhobracket}) give for any integer number $k \geq 2$
$$
\esp 
\left[
\left|
\overline{R}_{j}
-
\underline{R}_{j}
\right|^k
\right]
\leq 
\left(\frac{\overline{w}}{2}\right)^{k-2}
\esp 
\left[
\left(
\overline{R}_{j}
-
\underline{R}_{j}
\right)^2
\right]
\leq
\frac{k!}{8}
\left( \frac{\overline{w}}{2}\right)^{k-2}
\times
C
\frac{\left(1+t\right)^2}{ \underline{h}^{2(p+2)}}
\left(\eta^2+\eta+\eta^{1/2} \right).
$$
Hence 
(\ref{TBPbracks}) holds if $\eta$ satisfies
$$
\eta 
=
\frac{C}{3}
\min
\left(
\left(
\frac{\underline{h}^{2(p+2)}}{\left(1+t\right)^2}
\right)^{1/2}
\tau
,
\frac{\underline{h}^{2(p+2)}}{\left(1+t\right)^2}
\tau^2
,
\left(
\frac{\underline{h}^{2(p+2)}}{\left(1+t\right)^2}
\right)^{2}
\tau^4
\right).
$$
Recall now that $\tau < 1$, $t \geq 1$ and that  $\underline{h} \geq Cn^{-1/d}$ under Assumption \ref{K}. The bound (\ref{Neta}) for $N=\exp(H(\tau))$ gives taking $\eta$ as above
\begin{eqnarray*}
e^{H(\tau)}
& \leq &
\max
\left(
1,
\frac{C t^P}{
\min
\left(
\left(
\frac{\underline{h}^{2(p+2)}}{\left(1+t\right)^2}
\right)^{1/2}
\tau
,
\frac{\underline{h}^{2(p+2)}}{\left(1+t\right)^2}
\tau^2
,
\left(
\frac{\underline{h}^{2(p+2)}}{\left(1+t\right)^2}
\right)^{2}
\tau^4
\right)^{P+d+2}
}
\right)
\\
& \leq &
\max \left( 1, \frac{C t^3 n^{(4p+4)/d}}{\tau} \right)^{P+d+2}
.
\end{eqnarray*}
It then follows for $n$ large enough
$$
H(\tau)
\leq 
(P+d+2) \max \left( 0, \log\left(\frac{C t^3 n^{(4p+4)/d}}{\tau}\right) \right)
=
C
\left( 
   3 \log t 
   +
   \frac{4p+4}{d}\log n
    - 
   \log  \tau     
\right)
\leq 
C\log \left( \frac{tn}{\tau} \right),
$$
and (\ref{TBPH}) is proved. This ends the proof of the Lemma.
\eop

Let us now return to the proof of Proposition \ref{R1order}. Define $\mathbb{X} = (X_1,\cdots,X_n)$.
The definition of $\mathbb{R}_{n}^{1}$ and (\ref{R1}) give
\begin{eqnarray*}
\lefteqn{
\esp
\left[ 
\sup_{
  (\beta,\epsilon,\theta)
    \in
  \mathcal{B}(0,t_{\beta})
    \times
  \mathcal{B}(0,t_{\epsilon})
    \times
  \Theta^1
      }
\left|
  \mathbb{R}_{n}^{1}
      \left(
        \beta,\epsilon;\theta
      \right)
\right|
\right]
}
&&
\\
& = &
\esp
\left[
  \sup_{
    (\beta,\epsilon,\theta)
      \in
     \mathcal{B}(0,t_{\beta})
       \times 
      \mathcal{B}(0,t_{\epsilon})
       \times
      \Theta^1
       }
\left|
  \Sum_{i=1}^{n}
  \left(
    R_i
      \left(
        \beta, \epsilon;\theta
      \right)
    -
  \esp
    \left[
      R_i
        \left(
          \beta,\epsilon;\theta
        \right)
      |
      \mathbb{X}
    \right]
  \right)
\right| 
\right]
\\
& 
\leq
&
\esp
\left[ 
\sup_{(\beta,\epsilon,\theta)
      \in
     \mathcal{B}(0,t_{\beta}) 
       \times 
      \mathcal{B}(0,t_{\epsilon}) 
       \times
      \Theta^1
          }
\left|
  \Sum_{i=1}^{n}
    \left(
      R_i
        \left(
          \beta;\epsilon, \theta
        \right)
      -
      \esp
      \left[
        R_i
          \left(
            \beta;\epsilon,\theta
           \right)
       \right]
          \right)
\right|
\right]
\\
& &
+
\esp
\left[ 
\sup_{(\beta,\epsilon,\theta)
      \in
     \mathcal{B}(0,t_{\beta})
       \times 
      \mathcal{B}(0,t_{\epsilon}) 
       \times
       \Theta^1
          }
\left|
  \esp
    \left[
        \Sum_{i=1}^{n}
        \left(
          R_i
            \left(
              \beta,\epsilon;\theta
            \right)
          -
          \esp
          \left[
            R_i
              \left(
                \beta,\epsilon;\theta
              \right)
          \right]
        \right)
        \left| \mathbb{X} \right.
    \right]
\right|
\right]
\\
& 
\leq
& 
2\esp
\left[ 
\sup_{(\beta,\epsilon,\theta)
      \in
     \mathcal{B}(0,t_{\beta}) 
     \times 
     \mathcal{B}(0,t_{\epsilon})
     \times
     \Theta^1
           }
\left|
  \Sum^{n}_{i=1}
  \left(
    R_i
      \left(
        \beta,\epsilon;\theta
      \right)
    -
    \esp
    \left[
      R_i
        \left(
          \beta,\epsilon;\theta
        \right)
    \right]
      \right)
\right|
\right]
 .
\end{eqnarray*}
Let $H(\cdot)$, $\overline{\sigma}$ and $\overline{w}$ be as in Lemma \ref{Controlerho1}. Recall that $t_{\beta} + t_{\epsilon} \geq 1$ and that $\overline{\sigma} < 1 \leq n (t_{\beta} + t_{\epsilon})$ for $n$ large enough under the assumptions for $t_{\beta}$ and $t_{\epsilon}$ of the Proposition. It follows from Massart (2007, Theorem 6.8) that
$$
\esp
\left[ 
\sup_{(\beta,\epsilon,\theta)
      \in
     \mathcal{B}(0,t_{\beta}) 
     \times 
     \mathcal{B}(0,t_{\epsilon})
     \times
     \Theta^1
           }
\left|
  \Sum^{n}_{i=1}
  \!
  \left(
    R_i
      \left(
        \beta,\epsilon;\theta
      \right)
    -
    \esp
    \left[
      R_i
        \left(
          \beta,\epsilon;\theta
        \right)
    \right]
      \right)
\right|
\right]
\leq
C\left( n^{1/2} \!\! \int_{0}^{\overline{\sigma}} \!\!\! H(u)^{1/2}du + \left(\overline{w} + \overline{\sigma}\right)H\left(\overline{\sigma}\right)\right)
 .
$$
Since $\overline{\sigma}< 1$, Lemma \ref{Controlerho1} gives, for all $u$ in $(0,\overline{\sigma}]$, 
$H(u) \leq C\log (n(t_{\beta} + t_{\epsilon})/u)$. This gives
\begin{eqnarray*}
n^{1/2} \int_{0}^{\overline{\sigma}} H^{1/2} (u) du 
& \leq & 
(n \overline{\sigma} )^{1/2}
\left(
  \int_0^{\overline{\sigma}} H (u) du
\right)^{1/2}
\leq
C (n \overline{\sigma} )^{1/2}
\left(
  \int_0^{\overline{\sigma}} \log \left(\frac{n(t_{\beta} + t_{\epsilon})}{u}\right) du
\right)^{1/2}
\\
&  
 =
& 
C (n \overline{\sigma} )^{1/2}
\left(
\overline{\sigma}
\left(
\log \left(\frac{(t_{\beta} + t_{\epsilon})n}{\overline{\sigma}}\right)
+
1
\right)
\right)^{1/2}
\leq 
C n^{1/2} \overline{\sigma} \log^{1/2}\left( \frac{(t_{\beta} + t_{\epsilon})n}{\overline{\sigma}}\right)
 .
\end{eqnarray*}
The order for $\overline{\sigma}$ given in Lemma \ref{Controlerho1}, assumption on $t_{\beta} + t_{\epsilon}$ and Assumption \ref{K} give 
$$
\log\left(n(t_{\beta}+ t_{\epsilon})/\overline{\sigma}\right) 
\leq 
C
\log
\left(
\frac{n^{3/2} 
\left( n \underline{h}^d \right)^{1/4}
(t_{\beta}+ t_{\epsilon})^{1/2}}{t_{\epsilon}}
\right)
\leq
C
\log
\left(
\frac{n^{3/2} 
\left( n \underline{h}^d \right)^{1/2}
}{\log^{1/2} n}
\right)
\leq
C \log n.
$$    
Substituting gives
\begin{eqnarray*}
\lefteqn{
\esp
\left[ 
\sup_{
  (\beta,\epsilon,\theta)
    \in
  \mathcal{B}(0,t_{\beta})
    \times
  \mathcal{B}(0,t_{\epsilon})
    \times
  \Theta^1
      }
\left|
  \mathbb{R}_{n}^{1}
      \left(
        \beta,\epsilon;\theta
      \right)
\right|
\right]
\leq
C 
\left(
  n^{1/2} 
  \overline{\sigma} 
  \log^{1/2} 
  n
  +
  \left(
    \overline{\sigma} + \overline{w}
  \right)
  \log
  n
\right)
}
&&
\\
& \leq & 
C
  \frac{
     t_{\epsilon}(t_{\beta} + t_{\epsilon})^{1/2} 
              }
              {
               \left(n \underline{h}^d\right)^{1/4}
              }
\log^{1/2} n
\left(
  1 + 
     \log^{1/2} n
  \left(
    \frac{1}{n^{1/2}} 
      +
    \frac{(t_{\beta} + t_{\epsilon})^{1/2}}
    {t_{\epsilon}
    \left(
      n \underline{h}^d
    \right)^{1/4}}
  \right)
\right)
\leq  
 C
  \frac{
t_{\epsilon}(t_{\beta} 
+ t_{\epsilon})^{1/2}   
        }
        {
        \left(n \underline{h}^d\right)^{1/4}
        }
 \log^{1/2} 
        n
 .       
\eop
\end{eqnarray*}

\subsection{Proof of Proposition \protect{\ref{R2order}}}
The proof of Proposition \ref{R2order} follows the same steps of the proof of Proposition
\ref{R1order} and we only sketch it.
The integral expression of $R ( \beta,\epsilon;\theta)$ in
(\ref{Rho0.2}) and the expression (\ref{R2}) of ${\bf R}^2 ( \beta,\epsilon;\theta)$ give
$$
{\bf R}^2 (\beta,\epsilon;\theta)
= 
2
K_{h}
\left(
  X-x
\right)
\int_{\delta (\beta,\theta)}^{\delta (\beta,\theta) + \delta (\epsilon,\theta)}
\left(
       F
       \left(
              Q^* \left(X;\theta\right)
              +
              u
              |
              X
       \right)
-
       F
       \left(
              Q^* \left(X;\theta\right)
              |
              X
       \right)
\right)
du
-\frac{1}{2 n h^d}
\epsilon^T  {\bf J} (\theta) (\epsilon + 2\beta)
 .
$$
The definition (\ref{J}) of ${\bf J} (\theta)$ gives
\begin{eqnarray*}
\lefteqn{
{\bf R}^2(\beta,\epsilon;\theta)
}
&&
\\
& = &
2
K_{h}
\left(
  X-x
\right)
\int_{\delta (\beta,\theta)}^{\delta (\beta,\theta) + \delta (\epsilon,\theta)}
\left(
       F
       \left(
              Q^* \left(X;\theta\right)
              +
              u
              |
              X
       \right)
-
       F
       \left(
              Q^* \left(X;\theta\right)
              |
              X
       \right)
-
      u
      f
      \left(
              Q^* \left(X;\theta\right)
              |
              X
      \right)
\right)
du
\\
& = &
2
K_{h}
\left(
  X-x
\right)
\int_{\delta (\beta,\theta)}^{\delta (\beta,\theta) + \delta (\epsilon,\theta)}
u
\left\{
\int_0^1
\left(
       f
       \left(
              Q^* \left(X;\theta\right)
              +
              vu
              |
              X
       \right)
-
       f
       \left(
              Q^* \left(X;\theta\right)
              |
              X
       \right)
\right)
dv
\right\}
du
 .
\end{eqnarray*}
Define
$$
r(\beta;\theta)
=
2
K_{h}
\left(
  X-x
\right)
\int_{0}^{\delta (\beta,\theta)}
u
\left\{
\int_0^1
\left(
       f
       \left(
              Q^* \left(X;\theta\right)
              +
              vu
              |
              X
       \right)
-
       f
       \left(
              Q^* \left(X;\theta\right)
              |
              X
       \right)
\right)
dv
\right\}
du
$$
which is such that 
${\bf R}^2 (\beta,\epsilon;\theta)=r(\beta+\epsilon;\theta)-r(\beta;\theta)$. 
Since $| f(q+v|x)-f(q|x) |\leq L_0 |v|$ under Assumption \ref{F}, (\ref{Delta}) gives
\begin{eqnarray}
\left|
{\bf R}^2(\beta,\epsilon;\theta)
\right|
& \leq  &
K_h (X-x)
L_0\left|
\int_{\delta (\beta,\theta)}^{\delta (\beta,\theta) + \delta (\epsilon,\theta)}
u^2 du
\right|
\leq
C
K_{h}(X-x)
|\delta (\epsilon,\theta)|
(
  |\delta (\beta,\theta)| + |\delta (\epsilon,\theta)|
)^2
\nonumber \\
& \leq & 
C
\left\|{\bf U} \left(\frac{X-x}{h}\right)\right\|^3
K\left(\frac{X-x}{h}\right)
\frac{\left\|\epsilon\right\|\left(\left\|\beta\right\|+\left\|\epsilon\right\|\right)^2}{\left(nh^d\right)^{3/2}}
,
\label{MajorRho2}
\\
\left|
r(\beta;\theta)
\right|
& \leq &
CK_{h}(X-x)
|\delta (\beta,\theta)|^3
\leq
C
\left\|{\bf U} \left(\frac{X-x}{h}\right)\right\|^{3/2}
K\left(\frac{X-x}{h}\right)
\frac{\left\|\beta\right\|^3}{\left(nh^d\right)^{3/2}}
.
\nonumber
\end{eqnarray}
The latter inequality gives for all $\beta$ in $\mathcal{B} (0, t_{\beta}+t_{\epsilon})$ and all $\theta$ in $\Theta^1$
$$
\left|
r(\beta;\theta)
\right|
\leq
\frac{\overline{w}'}{2},
\quad
\overline{w}'
\asymp
\frac{(t_{\beta}+t_{\epsilon})^3}{\left(n\underline{h}^d\right)^{3/2}}.
$$
It follows from (\ref{MajorRho2}) that, for all $(\beta,\epsilon)$ in
$\mathcal{B} (0, t_{\beta}) \times \mathcal{B} (0, t_{\epsilon})$,
\begin{eqnarray*}
\lefteqn{
\Var
\left(
{\bf R}^2(\beta,\epsilon;\theta)
\right)
\leq
\esp 
\left[
{\bf R}^2(\beta,\epsilon;\theta)^2
\right]
}
&&
\\
& \leq &
C
\left(\frac{\left\|\epsilon\right\|\left(\left\|\beta\right\|+\left\|\epsilon\right\|\right)^2}{\left(nh^d\right)^{3/2}}
\right)^2
\int
\left\|{\bf U} \left(\frac{x'-x}{h}\right)\right\|^4
K^2\left(\frac{x'-x}{h}\right) 
f(x') dx'
\\
& \leq &
C
\frac{\left\|\epsilon\right\|^2\left(\left\|\beta\right\|+\left\|\epsilon\right\|\right)^4}{\left(nh^d\right)^{3}}
h^d
\int
\left\|{\bf U} \left(z\right)\right\|^4
K^2\left(z\right)
dz
\leq
\left(\overline{\sigma}'\right)^2
,
\quad
\overline{\sigma}'
\asymp
\frac{t_{\epsilon} \left(t_{\beta}+t_{\epsilon}\right)^2}{n^{3/2} \underline{h}^{d}}
.
\end{eqnarray*}
Then constructing brackets as in Lemma \ref{Controlerho1} and arguing as in the proof of Proposition \ref{R1order} give
\begin{eqnarray*}
\lefteqn{
\esp
\left[ 
\sup_{(\beta , \epsilon , \theta) 
  \in 
  \mathcal{B}(0,t_{\beta}) 
    \times 
  \mathcal{B}(0,t_{\epsilon}) 
    \times 
  \Theta^1} 
\left|
\mathbb{R}^2_n (\beta,\epsilon;\theta)
-
\esp
\left[
\mathbb{R}^2_n (\beta,\epsilon;\theta)
\right] 
\right|
\right]
}
& &
\\
&
\leq
& 
n^{1/2} 
\overline{\sigma}' 
\log^{1/2} 
  \left(
    \frac{ n(t_{\beta}+t_{\epsilon})}
    {\overline{\sigma}'}
  \right) 
  + 
  \left(
  \overline{\sigma}' + \overline{w}'
  \right) 
  \log
    \left(
      \frac{n(t_{\beta} + t_{\epsilon})}
      {\overline{\sigma}'}
    \right)
.
\end{eqnarray*}
Since (\ref{MajorRho2}) yields for all 
$(\beta , \epsilon , \theta) 
  \in 
  \mathcal{B}(0,t_{\beta}) 
    \times 
  \mathcal{B}(0,t_{\epsilon}) 
    \times 
  \Theta^1$
\begin{eqnarray*}
\left|
\esp
\left[
\mathbb{R}^2_n (\beta,\epsilon;\theta)
\right]
\right| 
& = &
\left|
n
\esp
\left[
\mathbf{R}^2 (\beta,\epsilon;\theta)
\right]
\right|
\leq
C
n
\esp
\left[
\left\|{\bf U} \left(\frac{X-x}{h}\right)\right\|^3
K\left(\frac{X-x}{h}\right)
\frac{\left\|\epsilon\right\|\left(\left\|\beta\right\|+\left\|\epsilon\right\|\right)^2}{\left(nh^d\right)^{3/2}}
\right]
\\
& \leq &
C
\frac{t_{\epsilon} \left(t_{\epsilon}+t_{\beta}\right)^2}{\left(n \underline{h}^{d}\right)^{1/2}}
, 
\end{eqnarray*}
substituting gives, using $t_{\beta} \geq 1$, $t_{\beta}/t_{\epsilon} =O \left( n \underline{h}^d/\log^{1/2} n\right)$ and Assumption \ref{K} which ensures $\log \left(n^{5/2} \underline{h}^d/t_{\epsilon}\right)= O (\log n)$,
\begin{eqnarray*}
\lefteqn{
\esp
\left[ 
\sup_{(\beta , \epsilon , \theta) 
  \in 
  \mathcal{B}(0,t_{\beta}) 
    \times 
  \mathcal{B}(0,t_{\epsilon}) 
    \times 
  \Theta^1} 
\left|
\mathbb{R}^2_n (\beta,\epsilon;\theta)
\right|
\right]
}
&&
\\
& \leq &
\esp
\left[ 
\sup_{(\beta , \epsilon , \theta) 
  \in 
  \mathcal{B}(0,t_{\beta}) 
    \times 
  \mathcal{B}(0,t_{\epsilon}) 
    \times 
  \Theta^1}
\left\{
\left|
\mathbb{R}^2_n (\beta,\epsilon;\theta)
-
\esp
\left[
\mathbb{R}^2_n (\beta,\epsilon;\theta)
\right] 
\right|
+
\esp
\left[
\mathbb{R}^2_n (\beta,\epsilon;\theta)
\right] 
\right\}
\right]
\\
& \leq &
C\frac{t_{\epsilon}\left(t_{\beta}+t_{\epsilon}\right)^2}{n \underline{h}^d}
\left(
1
+
\frac{t_{\beta}+t_{\epsilon}}{t_{\epsilon} \left(n \underline{h}^d\right)^{1/2}}
\right)
\log^{1/2}
\left(\frac{n^{5/2}\underline{h}^d}{t_{\epsilon}\left(t_{\beta}+t_{\epsilon}\right)}\right)
+
C
\frac{t_{\epsilon} \left(t_{\epsilon}+t_{\beta}\right)^2}{\left(n \underline{h}^{d}\right)^{1/2}}
\leq
\frac{t_{\epsilon} \left(t_{\epsilon}+t_{\beta}\right)^2}{\left(n \underline{h}^{d}\right)^{1/2}}
.
\eop
\end{eqnarray*}
\subsection{Proof of Lemma \protect{\ref{ControlJ}}}
Lemma \ref{Biais1} (iv) and Assumptions \ref{K} and \ref{F} give that there is a $C>0$ such that for all $\theta$ in $\Theta^1$ and all $i$,
\[
{\bf J}_i\left( \theta\right) \succ C {\bf M}_i\left(\theta\right)
 ,
\quad
{\bf M}_i\left(\theta\right) = 2K_h \left(X_i -x\right) {\bf U}\left(\frac{X_i - x}{h}\right) {\bf U}\left(\frac{X_i - x}{h}\right)^T .
\]
Hence for all $\theta$ in $\Theta^1$,
\begin{equation}
\frac{1}{nh^d} \Sum_{i=1}^{n} 
{\bf J}_i \left(\theta\right)
\succ
\frac{C}{nh^d}
\Sum_{i=1}^{n} {\bf M}_i
\left( \theta \right) = \mathbb{M}_n \left(\theta\right)
 .
\label{J2barM}
\end{equation}
The entries of $\mathbb{M}_n \left(\theta \right)$ write
\[
\frac{C}{nh^d}
\Sum_{i=1}^{n}
  \left(
    \frac{X_i-x}{h}
  \right)^{{\bf v}_1 +{\bf v}_2}
K\left(
  \frac{X_i -x}{h}
  \right)
 ,
0
  \leq
\left|{\bf v}_1\right|, \left|{\bf v}_2\right|
  \leq
p
 .
\]
Let $M(\theta)$ be the matrix with entries
\[
\frac{C}{h^d} \esp
          \left[
            \left(
              \frac{X -x}{h}
            \right)^{{\bf v}_1 + {\bf v}_2}
          K\left(
            \frac{X -x}{h}
           \right)
          \right]
,
0\leq |{\bf v}_1|, |{\bf v}_2| \leq p
.
\]
Arguing as in the proof of  Proposition \ref{R1order} for each of the entries of $\mathbb{M}_n \left(\theta \right)$ gives
\begin{equation*}
\sup_{\theta \in \Theta^1}
\|\mathbb{M}_n (\theta) - M(\theta) \|
=
o_{\Prob}\left(1\right)
.
\label{M2M}
\end{equation*}
Assumptions \ref{K}, \ref{F} and \ref{X} give, for  all ${\bf u}$ in $\Rit^P$, all $x$ in $\mathcal{X}_0$ and 
$\underline{h}$ small enough,
\begin{eqnarray*}
\lefteqn{
{\bf u}^T M(\theta) {\bf u}
= 
\frac{C}{h^d}
\esp
\Sum_{0\leq |{\bf v}_1|,|{\bf v}_2|\leq p}
u_{{\bf v}_1} u_{{\bf v}_2} 
\left(\frac{X-x}{h}\right)^{{\bf v}_1 + {\bf v}_2} K\left(\frac{X-x}{h}\right)
}
\\
& = &
C \Sum_{0\leq |{\bf v}_1|,|{\bf v}_2|\leq p}
u_{{\bf v}_1} u_{{\bf v}_2} \int z^{{\bf v}_1 + {\bf v}_2} K(z) f(x + h z)d z
= 
C\int
  \left(
    \Sum_{0\leq |{\bf v}| \leq p}
      u_{{\bf v}} z^{{\bf v}}
  \right)^2
  K(z)
  f(x + h z)
d z
\\
& \geq &
C \int_{\mathcal{B}(0,1)} \left( \sum_{0\leq |{\bf v}| \leq p} u_{\bf v} z^{\bf v} \right)^2 dz
\geq
C \left\|{\bf u}\right\|^2,
\end{eqnarray*}
where the last bound uses the fact that
$$
{\bf u}
\mapsto
\left(
\int_{\mathcal{B}(0,1)} \left( \sum_{0\leq |{\bf v}| \leq p} u_{\bf v} z^{\bf v} \right)^2 dz
\right)^{1/2}
$$
is a norm and that norms over $\Rit^P$ are equivalent.
Hence  (\ref{J2barM}) and $\|\mathbb{M}_n (\theta) - M(\theta) \|
     =
     o_{\Prob}\left(1\right)$ yield that there is a $\underline{\gamma} > 0$ such that
$
\inf_{\theta \in \Theta^1} 
\underline{\gamma}_n (\theta)
\geq
\inf_{\theta \in \Theta^1}  \inf_{\left\|{\bf u}\right\| = 1} {\bf u}^T \mathbb{M}_n (\theta){\bf u}
\geq \underline{\gamma} + o_{\prob} (1)
$ .
\eop
\subsection{Proof of Lemma \protect{\ref{Ordrebetan}}}
The first order condition (\ref{FOC}) implies that
$\esp [ {\bf S}_i (\theta) ] = 0$. 
Consider the ${\bf v}$ coordinate of ${\bf S}_i (\theta)$,
$$
{\bf S}_{{\bf v},i} (\theta)
=
2
\left\{
\ind
\left(
Y_i
\leq
Q^{*} (X_i;\theta)
\right)
-
\alpha
\right\}
\left(\frac{X_i-x}{h}\right)^{\bf v}
K
\left(\frac{X_i-x}{h}\right)
.
$$
Hence Assumptions \ref{K} and \ref{X} give, uniformly in $\theta \in \Theta^1$ and for all $i$,
\begin{eqnarray*}
\left|
\frac{{\bf S}_{{\bf v},i} (\theta)}{\left(nh^d\right)^{1/2}}
\right|
& \leq & 
\overline{w}'',
\quad
\overline{w}''
\asymp
\left(n\underline{h}^d\right)^{-1/2}
,
\\
\Var
\left(
\frac{{\bf S}_{{\bf v},i} (\theta)}{\left(n h^d\right)^{1/2}}
\right)
& \leq &
\esp
\left[\left( \frac{{\bf S}_{{\bf v},i} (\theta)}{\left(nh^d\right)^{1/2}} \right)^2\right]
\leq
\esp
\left[\left( \frac{\left(\frac{X_i-x}{h}\right)^{\bf v}
K
\left(\frac{X_i-x}{h}\right)}{\left(nh^d\right)^{1/2}} \right)^2\right]
=
\frac{h^d}{nh^d}
\int
\left( z^{\bf v}
K
\left(z\right) \right)^2
\\
& \leq &
\left( \overline{\sigma}'' \right)^2,
\quad
\overline{\sigma}''
\asymp
n^{-1/2}.
\end{eqnarray*} 
Hence arguing as in the proof of Proposition \ref{R1order} gives, under Assumption \ref{K},
\begin{equation*}
\esp
\left[
\sup_{\theta \in \Theta^1}
\left|
\frac{1}{(n\underline{h}^d)^{1/2}}
\sum_{i=1}^n
{\bf S}_{{\bf v},i} (\theta)
\right|
\right]
=
O
\left(
n^{1/2} \overline{\sigma}'' \log^{1/2}n
+
\left(\overline{\sigma}''+\overline{w}''\right)
\log^{1/2}n
\right)
=
O\left(\log^{1/2}n \right).
\end{equation*}
The Markov inequality then shows that the Lemma is proved.
\eop
}
\end{document}